\documentclass{birkjour}
\usepackage[noadjust]{cite}
\usepackage{amsthm, amssymb, amsmath} 
\usepackage{amssymb}
 \newtheorem{thm}{Theorem}[section]
 \newtheorem{cor}[thm]{Corollary}
 
 \newtheorem{prop}[thm]{Proposition}
 \theoremstyle{definition}
 \newtheorem{defn}[thm]{Definition}
 \theoremstyle{remark}
 \newtheorem{rem}[thm]{Remark}

 \newtheorem*{ex}{Example}
 \numberwithin{equation}{section}

\newcommand{\ed}{\end{document}}

\newcommand{\Spin}{\mathop{\mathrm{Spin}}}
\newcommand{\Pin}{\mathop{\mathrm{Pin}}}
\newcommand{\Cl}{\mathop{\mathrm{Cl}}}
\newcommand{\Aut}{\mathop{\mathrm{Aut}}}

\begin{document}
\bibliographystyle{spmpsci}
%
%
%
%
%
\submitted{December 26, 2020}
\received{December 26, 2020}
\revised{April 11, 2021}
%
%
%
%

\title[$H_4$ and the GA]
 {Clifford spinors and root system induction: $H_4$ and the Grand Antiprism}

\author[Dechant]{Pierre-Philippe Dechant}

\address{%
School of Science, Technology \& Health, York St John University, York YO31 7EX, United Kingdom\\{}\\
Department of Mathematics, University of York, Heslington YO10 5GE, United Kingdom\\{}\\
York Cross-disciplinary Centre for Systems Analysis, University of York, Heslington YO10 5GE, United Kingdom}

\email{ppd22@cantab.net}



\date{December 26, 2020}
\dedicatory{To John Horton Conway,  Richard Kenneth Guy and Michael Guy. }

\begin{abstract}
Recent work has shown that every 3D root system allows the construction of a correponding 4D root system via an `induction theorem'. 
In this paper, we look at the icosahedral case of $H_3\rightarrow H_4$ in detail and perform the calculations explicitly. 
Clifford algebra is used to perform group theoretic calculations based on the versor theorem and the Cartan-Dieudonn\'e theorem, giving a simple construction of the $\Pin$ and $\Spin$ covers.
Using this connection with $H_3$ via the induction theorem sheds light on geometric aspects of the $H_4$ root system (the $600$-cell) as well as other related polytopes and their symmetries, such as the famous Grand Antiprism and the snub 24-cell.
The uniform construction of root systems from 3D and the uniform procedure of splitting root systems with respect to subrootsystems into separate invariant sets allows further systematic insight into the underlying geometry.
All calculations are performed in the even subalgebra of $\Cl(3)$, including the construction of the Coxeter plane, which is used for visualising the complementary pairs of invariant polytopes, and are shared as supplementary computational work sheets. 
This approach therefore constitutes a more systematic and general way of performing calculations concerning groups, in particular reflection groups and root systems, in a Clifford algebraic framework. \end{abstract}

\maketitle

\section{Introduction}\label{sec_Intro}

Like many other mathematical concepts, symmetry groups have a rich structure including regular families and sporadic or exceptional phenomena \cite{he2015sporadic}. The exceptional symmetry group in 3D is icosahedral symmetry.  This fact has profound implications for the natural world around us, because objects consisting of identical building blocks that are `maximally symmetric'  display icosahedral symmetry. This includes most viruses and many fullerenes, as well as artificial nanocages in bionanotechnology and geodesic domes in architecture \cite{Dechant2021biochemist}. Even before any of these examples were known, icosahedral symmetry had inspired Plato to formulate a `unified theory of everything' in his dodecahedral `ordering principle of the universe'. This pattern of (exceptional) symmetries inspiring `grand unified theories' continues to this day, with $A_4=SU(5)$ in GUTs, and $E_8$ in string theory and GUTs, as well as $D_4=SO(8)$ and $B_4=SO(9)$ being critical in string and M theory. 

Traditionally, people seek to understand symmetries `top-down'. For instance $E_8$ includes $A_4$ and $H_4$, or $H_4$ includes $H_3$, so that people seek to understand the smaller groups as subgroups of the larger ones. In recent work \cite{Dechant2012CoxGA, Dechant2016Birth, Dechant2017e8, dechant2014SIGMAP} the author has shown that instead there is also a `bottom-up' view, by which e.g. $H_4$ and even $E_8$ can be constructed from $H_3$. In this sense the key to the larger exceptional symmetry groups is already contained in the smaller exceptional group. In particular, the author proved a uniform theorem that any 3D root system/reflection group induces a corresponding 4D root system/reflection group in a systematised way. The unusual abundance of exceptional symmetry structures in 4D could thus be based on (the accidentalness of) this construction, because it gives rise to the 4D exceptional objects $D_4$ (triality), $F_4$ (the largest 4D crystallographic group) and $H_4$ (the largest non-crystallographic group altogether). There is an immediate connection with Arnold's Trinities, mysterious connections between different triplets of exceptional objects throughout mathematics  \cite{Arnold1999symplectization, Arnold2000AMS,Dechant2018Trinity}.

There is therefore a lot of additional geometric insight to be gained from understanding this connection between 3D and 4D geometry, rather than looking at these phenomena from a 4D perspective alone. In fact, the link with 3D is much wider, including an ADE-type correspondence between 3D and 4D but also between 3D and ADE-type diagrams, in addition to the famous McKay correspondence \cite{Mckay1980graphs,Dechant2018Trinity}. The connection between 3D and 4D geometry arises because 3D reflections give rise to rotation groups via spinors. These 3D spinors themselves behave like 4D objects and can be shown to satisfy the root system axioms. This paper seeks to provide an example of concrete calculations performed entirely within the 3D Clifford algebra $\Cl(3)$ and its even subalgebra rather than in $\mathbb{R}^4$. There is of course also a connection with quaternions, because of an isomorphism. A lot of previous work has been done using quaternions and a purely algebraic description \cite{choi2018binary,Koca2009grand,Koca2006441,Koca2006H4,koca2007group,koca2012snub} with some somewhat haphazard results. However, we would argue that a lot of the deeper geometric insight, and its universality and systematic approach, have been lost by following an approach that is algebraically equivalent but is really not well suited for uncovering and understanding the underlying geometry and the generality of the situation. 

In this paper we therefore give a comprehensive and pedagogical exposition of how to perform group theoretic calculations for $H_3$ in a Clifford approach, leading to the detailed construction of  $H_3\rightarrow H_4$. We give explicit results of calculations in the paper, as well as making python Jupyter notebooks available as supplementary information for convenience and reproducibility. These detail the calculations, results and algorithms used based on the \verb|galgebra| python package \cite{Bromborsky2020} as well as some visualisations also from SageMath \cite{stein2008sage}. These notebooks are shared in the interest of open science, collegiality and reproducibility, and may be useful to the readers for adapting them for their own calculations. If they are useful please cite this paper along with the original software provided by e.g.  \cite{Bromborsky2020} and \cite{stein2008sage}.

After the basic construction of the $H_4$ root system (the $600$-cell) from $H_3$ we investigate various subrootsystems that arise within $H_4$. These groups and root systems are of course separately invariant within $H_4$, and the $120$ vertices of the $H_4$ root system can be split into two separately invariant sets by taking the complement of the subrootsystem in $H_4$. Conway and Guy found the Grand Antiprism computationally 55 years ago this year \cite{Conway1965Four}. Conway and Michael Guy's father, Richard Guy, himself a famous mathematician, populariser of mathematics and collaborator of Conway, have unfortunately passed away this year (2020). This paper is dedicated to their memory. The construction and the symmetries of the Grand Antiprism actually benefit from the construction from 3D as noted in \cite{Dechant2013Platonic}. In this paper we follow in detail how the $H_2\times H_2$ subgroup arises naturally within $H_4$ in the induction process. This subrootsystem is then used to split the $120$ vertices of $H_4$ into the set of $20$ roots of  $H_2\oplus H_2$ and the $100$ vertices of its complement, which is exactly the Grand Antiprism. A completely analogous construction works in a uniform way for other subrootsystems of $H_4$ that arise via the 3D construction, either as subrootsystems of $H_3$ or even subgroups of $2I$. The analogous cases include $D_4$ and the snub 24-cell, $A_1^4$, $A_2\oplus A_2$ and $A_4$.

We organise this paper as follows. We review some basics of Clifford algebras, reflection groups and root systems in Section \ref{sec_background}, leading to the Versor Theorem and the Induction Theorem. 
In Section \ref{sec_pin} we build on the Versor Theorem to set up a framework for explicit group theoretic calculations, including the construction of the $\Pin(H_3)$ and $\Spin(H_3)$ groups, and discussion of their conjugacy classes. 
Subrootsystems can arise either as even subgroups of the spinor group (here the binary icosahedral group $2I$) via subrootsystems of $H_3$ or generated via the inversion $e_1e_2e_3$, which is discussed in  Sections \ref{sec_group}
and  \ref{sec_rs}, respectively. 
The Coxeter plane is a convenient way to visualise root systems and other polytopes in any dimension, and the calculations in this case can also be entirely performed within the even subalgebra, as shown in Section \ref{sec_CoxPl}.
The next sections give detailed results for splitting the $H_4$ root system with respect to  various subrootsystems, yielding complementary pairs of invariant polytopes, including $H_2\oplus H_2$ and the Grand Antiprism in Section \ref{sec_GA}, $D_4$ and the snub 24-cell in Section \ref{sec_snub}, and examples from $A_1^4$, $A_2\oplus A_2$ and $A_4$ in Section \ref{sec_A2}.
In Section \ref{sec_concl} we conclude with a few words about the insights that can be gained from this novel approach via induction from 3D. 

\section{Root system induction}\label{sec_background}

The setting for reflection groups and root systems stipulates the existence of a vector space with an inner product, so without loss of generality one can construct the corresponding Clifford algebra. The Clifford algebraic framework used here is very standard, though we broadly follow  \cite{LasenbyDoran2003GeometricAlgebra, Hestenes1966STA}. We  define an algebra product via  the geometric product $xy=x\cdot y+x \wedge y$, where the inner product (given by the symmetric bilinear form) is the symmetric part $x\cdot y=\frac{1}{2}(xy+yx)$, and the wedge/exterior product the antisymmetric part $x\wedge y=\frac{1}{2}(xy-yx)$. These relations also mean that parallel vectors  commute whilst orthogonal vectors anticommute. The full $2^n$-dimensional algebra is generated via this geometric product, linearity and associativity. For our purposes we only consider the Clifford algebra of 3D $\Cl(3)$ generated by three orthogonal unit vectors $e_1$, $e_2$ and $e_3$, though some of the following statements hold under very general conditions. This yields an eight-dimensional vector space generated by the elements
$$
  \underbrace{\{1\}}_{\text{1 scalar}} \,\,\ \,\,\,\underbrace{\{e_1, e_2, e_3\}}_{\text{3 vectors}} \,\,\, \,\,\, \underbrace{\{e_1e_2=Ie_3, e_2e_3=Ie_1, e_3e_1=Ie_2\}}_{\text{3 bivectors}} \,\,\, \,\,\, \underbrace{\{I\equiv e_1e_2e_3\}}_{\text{1 trivector}},
$$
with an even subalgebra consisting of the scalar and bivectors, which is 4-dimensional. Note that for the orthogonal unit vectors e.g. $e_1e_2=e_1\wedge e_2$. We follow the \verb|galgebra| \LaTeX  { }output, which has the wedge version.

Root systems and in particular their simple roots are convenient objects to characterise reflection and Coxeter groups. We therefore briefly introduce the relevant terminology here:
\begin{defn}[Root system] \label{DefRootSys}
A \emph{root system} is a collection $\Phi$ of non-zero (root)  vectors $\alpha$ that span an $n$-dimensional Euclidean vector space $V$ endowed with a positive definite bilinear form, that satisfies the  two axioms:
\begin{enumerate}
\item $\Phi$ only contains a root $\alpha$ and its negative, but no other scalar multiples: $\Phi \cap \mathbb{R}\alpha=\{-\alpha, \alpha\}\,\,\,\,\forall\,\, \alpha \in \Phi$. 
\item $\Phi$ is invariant under all reflections corresponding to root vectors in $\Phi$: $s_\alpha\Phi=\Phi \,\,\,\forall\,\, \alpha\in\Phi$. 
The reflection $s_\alpha$ in the hyperplane with normal vector $\alpha$ is given by $$s_\alpha: x\rightarrow s_\alpha(x)=x - 2\frac{(x\cdot\alpha)}{(\alpha\cdot \alpha)}\alpha,$$ where $(x \cdot y)$ denotes the inner product on $V$.
\end{enumerate}
\end{defn}
Unlike other popular conventions here we assume unit normalisation for all our root vectors for later convenience. 

\begin{prop}[Reflection]\label{prop_refl}
In Clifford algebra, the reflection formula simplifies to $$s_\alpha: x\rightarrow s_\alpha(x)=-\alpha x \alpha$$ for normalised root vectors, so that double-sided (`sandwich') application of a root vector generates the corresponding reflection. 
\end{prop}

\begin{proof}
Using the Clifford form for the inner product $x\cdot y=\frac{1}{2}(xy+yx)$ in the (simple) reflection formula in Definition \ref{DefRootSys}   $s_i: x\rightarrow s_i(x)=x - 2\frac{(x\cdot\alpha)}{(\alpha\cdot \alpha)}\alpha$ and assuming unit normalisation of roots $\alpha_i\cdot \alpha_i=1$ yields the much simplified version 
$$s_i: x\rightarrow s_i(x)=x-2\cdot \frac{1}{2}(x \alpha_i+ \alpha_ix) \alpha_i=x-x \alpha_i^2- \alpha_ix \alpha_i=-\alpha_ix\alpha_i.$$
\end{proof}

\begin{prop}[Double cover]
In Clifford algebra, the reflections $\alpha$ and $-\alpha$ doubly cover the same reflection $s_\alpha$.\label{prop_double_cover}
\end{prop}
\begin{proof}
Straightforward, since due to the double-sided application the signs cancel out. 
\end{proof}

A subset $\Delta$ of $\Phi$, called \emph{simple roots} $\alpha_1, \dots, \alpha_n$ is sufficient to express every element of $\Phi$ via linear combinations with coefficients of the same sign. 
For a \emph{crystallographic} root system, these are $\mathbb{Z}$-linear combinations, whilst for the \emph{non-crystallographic} root systems one needs to consider certain extended integer rings. For instance  for $H_2$, $H_3$ and $H_4$ one has the extended integer ring $\mathbb{Z}[\tau]=\lbrace a+\tau b| a,b \in \mathbb{Z}\rbrace$, where $\tau$ is   the golden ratio $\tau=\frac{1}{2}(1+\sqrt{5})=2\cos{\frac{\pi}{5}}$, and $\sigma$ is its Galois conjugate $\sigma=\frac{1}{2}(1-\sqrt{5})$ (the two solutions to the quadratic equation $x^2=x+1$), and linear combinations are with respect to this $\mathbb{Z}[\tau]$. This integrality property of the crystallographic root systems  (types $A$-$G$) leads to an associated lattice which acts as a root lattice for Lie algebras, which are named accordingly. In contrast, no such lattice exists for the non-crystallographic groups (types $H$ and $I$), which accordingly do not have associated Lie algebras, and are perhaps less familiar as a result.

	The reflections corresponding to simple roots are also called simple reflections.  The geometric structure of the set of simple roots encodes the properties of the reflection group and is summarised in the Cartan matrix and Coxeter-Dynkin diagrams, which contain the geometrically invariant information of the root system as follows: 
\begin{defn}[Cartan matrix and Coxeter-Dynkin diagram] The  \emph{Cartan matrix} of a set of simple roots $\alpha_i\in\Delta$ is defined as the matrix
\begin{equation}\label{CM}
	A_{ij}=2\frac{(\alpha_i\cdot \alpha_j)}{(\alpha_i\cdot \alpha_i)}.
\end{equation}
A graphical representation of the geometric content is given by \emph{Coxeter-Dynkin diagrams}, in which nodes correspond to simple roots, orthogonal roots are not connected, roots at $\frac{\pi}{3}$ have a simple link, and other angles $\frac{\pi}{m}$ have a link with a label $m$. 
\end{defn}

\begin{ex}
The Cartan matrices for $H_3$ and $H_4$ are respectively given by
$$A \left(H_3\right) = \begin{pmatrix}
	    2&-1&0
	\\ -1&2&-\tau
	\\  0&-\tau&2
	
	 \end{pmatrix}, A \left(H_4\right) = \begin{pmatrix}
	    2&-1&0&0
	\\ -1&2&-1&0 
	\\  0&-1&2&-\tau
	\\  0&0&-\tau&2
	 \end{pmatrix}.$$
	 Possible choices of simple roots are e.g. $$a_1 =e_2, a_2 = \frac{1}{2}(-\tau e_1-e_2-(\tau-1)e_3), a_3=e_1$$ and $$\alpha_1=\frac{1}{2}(\tau e_1-e_2+(\tau -1) e_4), \,\,\,\alpha_2=e_2,\,\,\, \alpha_3=-\frac{1}{{2}}((\tau -1)e_1+e_2+\tau e_3) \text{ and } \alpha_4=e_3.$$
	 \end{ex}
	 
\begin{rem}
What is a slight drawback of the diagrammatic approach is that it is rather unobvious which subgroups are contained e.g. in $H_4$ -- at least exhaustively. We will see some non-obvious examples later on. \end{rem}

Root systems and simple roots are therefore convenient paradigms for considering reflection groups: each root vector defines a hyperplane that it is normal to and thereby acts as a generator of a reflection in that hyperplane. In Clifford algebra this root vector in fact directly acts as a reflection generator via the geometric product. Multiplying together such simple reflections 
$s_i: x\rightarrow s_i(x)=x - 2\frac{(x|\alpha_i)}{(\alpha_i|\alpha_i)}\alpha_i=-\alpha_ix\alpha_i$
therefore generates a reflection group. This is in fact a Coxeter group, since the simple reflections $s_i$ satisfy the defining relations:
\begin{defn}[Coxeter group] A \emph{Coxeter group} is a group generated by a set of involutory generators $s_i, s_j \in S$ subject to relations of the form $(s_is_j)^{m_{ij}}=1$ with $m_{ij}=m_{ji}\ge 2$ for $i\ne j$. 
\end{defn}
\begin{defn}[Coxeter element and Coxeter number] The product of all the simple reflections in some order is called a Coxeter element. All such elements are conjugate and as such their order is well-defined and called the Coxeter number. 
\end{defn}

These reflection groups are built up in Clifford algebra by performing successive multiplication with the unit vectors defining the reflection hyperplanes via `sandwiching'
\begin{equation}
s_1\dots s_k: x\rightarrow s_1\dots s_k(x)=(-1)^k\alpha_1\dots\alpha_kx\alpha_k\dots\alpha_1=:(-1)^k Ax\tilde{A},\label{GACox}
\end{equation}
where the tilde denotes the \emph{reversal} of the order of the constituent vectors in the product $A=\alpha_1\dots\alpha_k$. In order to study the groups of reflections one therefore only needs to consider  products of root vectors in the Clifford algebra, which form a multivector group under the geometric product and yield a Pin double cover of the corresponding reflection group \cite{Porteous1995Clifford}. The inverse of each group element is of course simply given by the reversal, because of the assumed normalisation condition. Since $\alpha_i$ and $-\alpha_i$ encode the same reflection, products of unit vectors are double covers of the respective orthogonal transformation, as $A$ and $-A$ encode the same transformation. We call even products $R$, i.e. products of an even number of vectors, `spinors' or `rotors', and a general product $A$ `versors' or `pinors'. They form the $\Pin$ group and constitute a double cover of the orthogonal group, whilst the even products form the double cover of the special orthogonal group, called the $\Spin$ group.
Clifford algebra therefore provides a particularly natural and simple construction of the $\Spin$ groups.

In fact, more general groups can be constructed in this way because of the fundamental importance of reflections according to the Cartan-Dieudonn\'e Theorem \cite{Garling2011Clifford}.
\begin{thm}[Cartan-Dieudonn\'e Theorem]
Every orthogonal transformation in an $n$-dimensional symmetric bilinear space can be described as the composition of at most $n$ reflections.\label{thrm_CD}
\end{thm}
The above approach to group theory via multivector groups \eqref{GACox} is therefore a much more general way of doing group theory.

In Clifford algebra, instead of using matrices to perform linear transformations one can use spinors/rotors/pinors/versors to perform linear transformations that leave the inner product invariant i.e. orthogonal transformations \cite{Hestenes2002CrystGroups,Hestenes2002PointGroups,Hitzer2010CLUCalc}. Here the normalisation condition has been dropped as long as the vectors are non-null 
since the  inverse of multiplication with a non-null vector $x$ is simply $x^{-1}=\frac{x}{|x|^2}$ since $xx=x\cdot x = |x|^2$ (in the positive signature spaces we will consider there are no null vectors anyway). Therefore the multivector $A$ that is a product of vectors is invertible and preserves the inner product, though the inverse is no longer just given by the reverse:
\begin{thm}[Versor Theorem]
Every orthogonal transformation $\underline{A}$ can be expressed in the canonical form
$\underline{A}:x \rightarrow x'=\underline{A}(x) =\pm A^{-1}xA$ where $A$ is a versor and the sign is its parity. \label{thrm_versor}
\end{thm}

A concept familiar from abstract group theory via generators and relations is that group elements can be written as words in the generators. It is noteworthy here that in contrast to for instance the Coxeter group elements as words in the generators $s_\alpha$, in the Clifford algebra approach the root vectors are \emph{directly} generators for the Pin double cover under multiplication with the geometric product. To stress this slight distinction we call these `generator paths' for each versor $v$, with $\pm v$ each versor corresponding to one Coxeter group word.  

Therefore Clifford algebras and root systems are frameworks that perfectly complement each other since performing reflections in Clifford algebras is so simple and only assumes the structure of a vector space with an inner product that is already given in the root system definition. Therefore Clifford algebras are perhaps the most natural framework for studying reflection groups and root systems, and through the above arguments also more general groups \cite{Dechant2015ICCA}. In the next section we perform a detailed computation of the $\Pin$ and $\Spin$ covers of the icosahedral groups to illustrate the principles. Following the above two frameworks, given (simple) roots in a root system one can start multiplying these together using the geometric product. General products will be in $\Pin$ whilst even products are in $\Spin$. For now we will concentrate on the $\Spin$ and the even subalgebra.

\begin{prop}[$O(4)$-structure of spinors]\label{HGA_O4}
The space of $\Cl(3)$-spinors has a 4D Euclidean structure.
\end{prop}
\begin{proof}
	For a spinor $R=a_0+a_1e_2e_3+a_2e_3e_1+a_3e_1e_2$, the norm is given by $R\tilde{R}=a_0^2+a_1^2+a_2^2+a_3^2$, and the inner product between two spinors $R_1$ and $R_2$ is $(R_1, R_2)=\frac{1}{2}(R_1\tilde{R}_2+R_2\tilde{R}_1)$.
\end{proof}

\begin{rem}
For rotors, the inner product $\frac{1}{2}(R_1\tilde{R}_2+R_2\tilde{R}_1)$ is of course invariant under $R_i\rightarrow RR_i\tilde{R}$ but also just under $R_i\rightarrow RR_i$.\end{rem}

From the double cover property of Proposition \ref{prop_double_cover} we have the following corollary: 
\begin{cor}[Discrete spinor groups]\label{HGA_disspin}
	Discrete spinor groups are of even order since if a spinor $R$ is contained in the group then so is $-R$ since it encodes the same orthogonal transformation.
\end{cor}

\begin{prop}[Spin group closure properties]
Spin groups are closed under: \label{prop_spin_closure}
\begin{itemize}
\item multiplication using the geometric product
\item reversal
\item multiplication by $-1$
\end{itemize}
\end{prop}

\begin{proof}
Straightforward:
\begin{itemize}
\item by definition of multivector groups via the geometric product
\item by the inverse element group axiom since reversal is equivalent to the inverse
\item by Corollary \ref{HGA_disspin} both $R$ and $-R$ are contained in the group 
\end{itemize}
\end{proof}	

Following the formula for fundamental reflections from Definition \ref{DefRootSys} one can likewise define reflections on this spinor space with respect to the inner product between spinors. 
\begin{prop}[Spin reflections]
Reflections between spinors using the spinor inner product are given by $$R_2\rightarrow R_2'=  -R_1\tilde{R}_2R_1.$$\label{prop_spin_refl}
\end{prop}

\begin{proof}
 In analogy to Proposition \ref{prop_refl}, for normalised spinors $R_1$ and $R_2$ and using the definition of the spinor inner product from Proposition \ref{HGA_O4} this amounts to 
$R_2\rightarrow R_2'=R_2-2(R_1, R_2)/(R_1, {R}_1) R_1 =R_2-(R_1\tilde{R}_2+R_2\tilde{R}_1) R_1=R_2-R_1\tilde{R}_2 R_1-R_2\tilde{R}_1 R_1=  -R_1\tilde{R}_2R_1$.
\end{proof}

\begin{prop}[3D spinor -- 4D vector correspondence]
Spinor reflections in the spinor $R=a_0+a_1e_2\wedge e_3- a_2e_1\wedge e_3 + a_3e_1\wedge e_2$ are equivalent to 4D reflections in the 4D vector $(a_0, a_3, -a_2, a_1)$.
\end{prop}

\begin{proof}
By direct calculation or see supplementary material. 
\end{proof}
In practice, the exact mapping of the components is a matter of convention and often irrelevant, since most root systems contain roots up to $\pm$ and (cyclic) permutations anyway. 

\begin{thm}[Induction Theorem]\label{thrm_ind}
	Any rank-3 root system induces a  root system of rank 4.
\end{thm}

\begin{proof}
A root system in three dimensions $\Phi^{(3)}$ gives rise to a group of spinors by taking even products of the root vectors. 
From Corollary \ref{HGA_disspin}, this group contains $-R$ if it contains $R$, and therefore satisfies the first root system axiom from Definition \ref{DefRootSys}. The set of spinors has a 4D Euclidean structure by Proposition \ref{HGA_O4} and can thus be treated as a collection of 4D vectors $\Phi^{(4)}$ with the inner product as given in the Proposition. It remains to show that this collection of vectors $\Phi^{(4)}$ is invariant under reflections (axiom 2), which is satisfied by  Proposition \ref{prop_spin_refl}. 
\end{proof}
Closure of the root system is thus ensured by closure of the spinor group. This also has very interesting consequences for the automorphism group of these spinorial root systems, which contains two factors of the spinor group acting from the left and the right \cite{Dechant2013Platonic} (in this sense, the above closure under reflections amounts to a certain twisted conjugation). 

There is a limited number of cases which we can just enumerate. The 3D root systems are listed in Table \ref{RSenumeration} along with the 4D root systems that they induce as well as the intermediate spinor groups (the binary polyhedral groups). In this article, we focus on the case $H_3 \rightarrow H_4$.

\begin{table}
	\caption{List of the rank 3 root systems and  their induced root systems in four dimensions, as well as the binary polyhedral groups that act as the spinor group intermediaries.}
\begin{centering}\begin{tabular}{|clc|c}
\hline 
start root system&induced root system& binary polyhedral group
\tabularnewline
\hline 
$A_1^3$&$A_1^4$& $Q$
\tabularnewline
$A_1\oplus I_2(n)$&$I_2(n)\oplus I_2(n)$ & $<2,2,n>$
\tabularnewline
$A_3$&$D_4$ &$2T$
\tabularnewline
$B_3$&$F_4$ & $2O$
\tabularnewline
$H_3$&$H_4$ & $2I$
\tabularnewline
\hline 
\end{tabular}\end{centering}\label{RSenumeration}
\end{table}
\begin{defn}[Subrootsystem] By a subrootsystem $\Phi_1$ of a root system $\Phi_2$ we mean a subset $\Phi_1$ of the collection of vectors $\Phi_2$ that itself satisfies the root system axioms. 
\end{defn}
From the Induction Theorem \ref{thrm_ind} we immediately get:
\begin{cor}[Induced subrootsystems]\label{ind_subRS}
	A subrootsystem $\Phi_1^{(3)}$ of a root system $\Phi_2^{(3)}$ induces a subrootsystem $\Phi_1^{(4)}$ of the induced root system $\Phi_2^{(4)}$. 
\end{cor}

Any subrootsystem of $H_3$ therefore induces a subrootsystem of $H_4$. For instance, the $A_1^3$ inside $H_3$ induces the rather boring $A_1^4$ in $H_4$. Similarly, $A_2$ and $H_2$ are contained, if rather boring as 2D root systems. If $H_3$  contained $A_1\oplus A_2$ and $A_1\oplus H_2$ subrootsytems (which it doesn't), then this would lead to the doubling $A_2\oplus A_2$ and $H_2\oplus H_2$ inside the $H_4$. This is not quite the case, but nearly so, which we will return to later. Similarly, it follows from the Induction Theorem that any even subgroup of a spinor group will also yield a subrootsystem. We will explore these points in the  Sections \ref{sec_group}
and  \ref{sec_rs} which follow the next section, where we will discuss the multivector group calculation framework.

\section{Pin Group and Spin Group}\label{sec_pin}

In this section we pick up on the practical implications of the versor theorem \ref{thrm_versor} and the Cartan-Dieudonn\'e theorem in order to do explicit calculations in group theory 
in the Clifford algebra approach. We use the concrete example of the $H_3$ reflection/Coxeter group. The icosahedral rotation group (in $SO(3)$) is the alternating group $A_5$ of order 60, also known simply as $I$ (which we will avoid due to the pseudoscalar often being denoted by that too). This group is of course doubly covered in $\Spin(3)$  by its spin double cover,  with its nice  Clifford algebra construction via the reflection formula in Proposition \ref{prop_refl}. We might denote this group by $\Spin(H_3)$ but it is also commonly known as the binary icosahedral group $2I$. From the induction theorem \ref{thrm_ind} of the previous section the elements of this group of course give the $120$ roots of the $H_4$ root system. The rotational group is of course also doubly covered in $O(3)$ by its double cover $H_3=A_5 \times \mathbb{Z}_2$. Both double covers are of course of order $120$. $H_3$ however is itself also doubly covered in  $\Pin(3)$ by a group of order $240$, which doesn't have a common name since it is simply $H_3\times \mathbb{Z}_2=A_5\times \mathbb{Z}_2\times \mathbb{Z}_2$ but which we might for consistency call $\Pin(H_3)$. 

For the reader's convenience and for reproducibility in Tables \ref{tab:CC_spin1} to \ref{tab:CC_pin2}  we list  the different group elements of $\Pin(H_3)$ in our Clifford approach explicitly. The python Jupyter notebooks in the supplementary material contain the algorithms used, which are based on the \verb|galgebra| software package \cite{Bromborsky2020}. For convenience we group the elements as whole conjugacy classes but we give the order in which the elements are generated by repeated application of the generating simple roots. This gives a number for each group element as a reference for ease of access, as well as the word in the generators (the `generator path') that generates this particular group element in terms of the $H_3$ simple roots/generators $a_1, a_2, a_3$. For the simple roots of $H_3$ we pick
$$a_1 =e_2, a_2 = \frac{1}{2}(-\tau e_1-e_2-(\tau-1)e_3), a_3=e_1.$$ Note that the wedge could be omitted since we have picked the orthogonal unit vectors $e_1, e_2, e_3$ so the wedge product is synonymous with the full geometric product. We multiply all group elements by $2$ to save clutter.

\tiny

\begin{table}
\begin{tabular}{|c|c|c|c|c|c|c}
   \hline
   \textbf{Order}  & \textbf{Number} & \textbf{Element} $\times 2$ & \textbf{Generator path}  \\
  \hline 

 $ 1 $ &  $ 4 $ & $ 2  $ & 11 \\
  \hline 

 $ 2 $ &  $ 26 $ & $ -2  $ & 1313 \\
  
  \hline

  \hline
 $ 3 $ &  $ 5 $ & $ -1  + \tau  e_{1}\wedge  e_{2} + \sigma  e_{2}\wedge  e_{3}  $ & 12 \\ 
$ 3 $ &  $ 7 $ & $ -1  -\tau  e_{1}\wedge  e_{2} -\sigma  e_{2}\wedge  e_{3}  $ & 21 \\ 
$ 3 $ &  $ 36 $ & $ -1  -\tau  e_{1}\wedge  e_{2} + \sigma  e_{2}\wedge  e_{3}  $ & 3123 \\ 
$ 3 $ &  $ 38 $ & $ -1  + \tau  e_{1}\wedge  e_{2} -\sigma  e_{2}\wedge  e_{3}  $ & 3213 \\ 
$ 3 $ &  $ 82 $ & $ -1  + e_{1}\wedge  e_{2} + e_{1}\wedge  e_{3} -e_{2}\wedge  e_{3}  $ & 231232 \\ 
$ 3 $ &  $ 83 $ & $ -1  -e_{1}\wedge  e_{2} -e_{1}\wedge  e_{3} + e_{2}\wedge  e_{3}  $ & 232132 \\ 
$ 3 $ &  $ 124 $ & $ -1  -e_{1}\wedge  e_{2} + e_{1}\wedge  e_{3} + e_{2}\wedge  e_{3}  $ & 12312321 \\ 
$ 3 $ &  $ 125 $ & $ -1  + e_{1}\wedge  e_{2} -e_{1}\wedge  e_{3} -e_{2}\wedge  e_{3}  $ & 12321321 \\ 
$ 3 $ &  $ 128 $ & $ -1  -e_{1}\wedge  e_{2} + e_{1}\wedge  e_{3} -e_{2}\wedge  e_{3}  $ & 12323123 \\ 
$ 3 $ &  $ 131 $ & $ -1  -e_{1}\wedge  e_{2} -e_{1}\wedge  e_{3} -e_{2}\wedge  e_{3}  $ & 13213232 \\ 
$ 3 $ &  $ 134 $ & $ -1  -\sigma  e_{1}\wedge  e_{2} + \tau  e_{1}\wedge  e_{3}  $ & 21321323 \\ 
$ 3 $ &  $ 137 $ & $ -1  + e_{1}\wedge  e_{2} + e_{1}\wedge  e_{3} + e_{2}\wedge  e_{3}  $ & 23213213 \\ 
$ 3 $ &  $ 141 $ & $ -1  + \sigma  e_{1}\wedge  e_{2} -\tau  e_{1}\wedge  e_{3}  $ & 32132132 \\ 
$ 3 $ &  $ 143 $ & $ -1  + e_{1}\wedge  e_{2} -e_{1}\wedge  e_{3} + e_{2}\wedge  e_{3}  $ & 32132321 \\ 
$ 3 $ &  $ 170 $ & $ -1  + \sigma  e_{1}\wedge  e_{2} + \tau  e_{1}\wedge  e_{3}  $ & 1213213231 \\ 
$ 3 $ &  $ 172 $ & $ -1  -\sigma  e_{1}\wedge  e_{3} -\tau  e_{2}\wedge  e_{3}  $ & 1213231232 \\ 
$ 3 $ &  $ 174 $ & $ -1  + \sigma  e_{1}\wedge  e_{3} -\tau  e_{2}\wedge  e_{3}  $ & 1232132312 \\ 
$ 3 $ &  $ 177 $ & $ -1  -\sigma  e_{1}\wedge  e_{2} -\tau  e_{1}\wedge  e_{3}  $ & 1321321321 \\ 
$ 3 $ &  $ 184 $ & $ -1  -\sigma  e_{1}\wedge  e_{3} + \tau  e_{2}\wedge  e_{3}  $ & 2132312321 \\ 
$ 3 $ &  $ 187 $ & $ -1  + \sigma  e_{1}\wedge  e_{3} + \tau  e_{2}\wedge  e_{3}  $ & 2321323121 \\ 
  \hline

      \hline 
 $ 4 $ &  $ 6 $ & $ -2  e_{1}\wedge  e_{2}  $ & 13 \\
$ 4 $ &  $ 9 $ & $ 2  e_{1}\wedge  e_{2}  $ & 31 \\
$ 4 $ &  $ 30 $ & $ \tau  e_{1}\wedge  e_{2} -\sigma  e_{1}\wedge  e_{3} -  e_{2}\wedge  e_{3}  $ & 2132 \\
$ 4 $ &  $ 31 $ & $ -\tau  e_{1}\wedge  e_{2} + \sigma  e_{1}\wedge  e_{3} +   e_{2}\wedge  e_{3}  $ & 2312 \\
$ 4 $ &  $ 66 $ & $ -\tau  e_{1}\wedge  e_{2} -\sigma  e_{1}\wedge  e_{3} +   e_{2}\wedge  e_{3}  $ & 121321 \\
$ 4 $ &  $ 68 $ & $ \tau  e_{1}\wedge  e_{2} + \sigma  e_{1}\wedge  e_{3} -  e_{2}\wedge  e_{3}  $ & 123121 \\
$ 4 $ &  $ 89 $ & $ -\tau  e_{1}\wedge  e_{2} + \sigma  e_{1}\wedge  e_{3} -  e_{2}\wedge  e_{3}  $ & 321323 \\
$ 4 $ &  $ 91 $ & $ \tau  e_{1}\wedge  e_{2} -\sigma  e_{1}\wedge  e_{3} +   e_{2}\wedge  e_{3}  $ & 323123 \\
$ 4 $ &  $ 129 $ & $ -\tau  e_{1}\wedge  e_{2} -\sigma  e_{1}\wedge  e_{3} -  e_{2}\wedge  e_{3}  $ & 13213213 \\
$ 4 $ &  $ 130 $ & $ \tau  e_{1}\wedge  e_{2} + \sigma  e_{1}\wedge  e_{3} +   e_{2}\wedge  e_{3}  $ & 13213231 \\
$ 4 $ &  $ 139 $ & $   e_{1}\wedge  e_{2} + \tau  e_{1}\wedge  e_{3} + \sigma  e_{2}\wedge  e_{3}  $ & 23213232 \\
$ 4 $ &  $ 140 $ & $ -  e_{1}\wedge  e_{2} -\tau  e_{1}\wedge  e_{3} -\sigma  e_{2}\wedge  e_{3}  $ & 23231232 \\
$ 4 $ &  $ 175 $ & $ -  e_{1}\wedge  e_{2} + \tau  e_{1}\wedge  e_{3} -\sigma  e_{2}\wedge  e_{3}  $ & 1232132321 \\
$ 4 $ &  $ 176 $ & $   e_{1}\wedge  e_{2} -\tau  e_{1}\wedge  e_{3} + \sigma  e_{2}\wedge  e_{3}  $ & 1232312321 \\
$ 4 $ &  $ 181 $ & $ -\sigma  e_{1}\wedge  e_{2} +   e_{1}\wedge  e_{3} -\tau  e_{2}\wedge  e_{3}  $ & 2132132132 \\
$ 4 $ &  $ 182 $ & $ \sigma  e_{1}\wedge  e_{2} -  e_{1}\wedge  e_{3} + \tau  e_{2}\wedge  e_{3}  $ & 2132132312 \\
$ 4 $ &  $ 186 $ & $ -  e_{1}\wedge  e_{2} -\tau  e_{1}\wedge  e_{3} + \sigma  e_{2}\wedge  e_{3}  $ & 2321321323 \\
$ 4 $ &  $ 188 $ & $   e_{1}\wedge  e_{2} + \tau  e_{1}\wedge  e_{3} -\sigma  e_{2}\wedge  e_{3}  $ & 2321323123 \\
$ 4 $ &  $ 209 $ & $ \sigma  e_{1}\wedge  e_{2} +   e_{1}\wedge  e_{3} + \tau  e_{2}\wedge  e_{3}  $ & 121321321321 \\
$ 4 $ &  $ 211 $ & $ -\sigma  e_{1}\wedge  e_{2} -  e_{1}\wedge  e_{3} -\tau  e_{2}\wedge  e_{3}  $ & 121321323121 \\
$ 4 $ &  $ 213 $ & $ -  e_{1}\wedge  e_{2} + \tau  e_{1}\wedge  e_{3} + \sigma  e_{2}\wedge  e_{3}  $ & 123213213213 \\
$ 4 $ &  $ 214 $ & $   e_{1}\wedge  e_{2} -\tau  e_{1}\wedge  e_{3} -\sigma  e_{2}\wedge  e_{3}  $ & 123213213231 \\
$ 4 $ &  $ 221 $ & $ \sigma  e_{1}\wedge  e_{2} -  e_{1}\wedge  e_{3} -\tau  e_{2}\wedge  e_{3}  $ & 321321321323 \\
$ 4 $ &  $ 222 $ & $ -\sigma  e_{1}\wedge  e_{2} +   e_{1}\wedge  e_{3} + \tau  e_{2}\wedge  e_{3}  $ & 321321323123 \\
$ 4 $ &  $ 233 $ & $ -2  e_{2}\wedge  e_{3}  $ & 12132132132132 \\
$ 4 $ &  $ 234 $ & $ 2  e_{2}\wedge  e_{3}  $ & 12132132132312 \\
$ 4 $ &  $ 235 $ & $ -\sigma  e_{1}\wedge  e_{2} -  e_{1}\wedge  e_{3} + \tau  e_{2}\wedge  e_{3}  $ & 12321321321323 \\
$ 4 $ &  $ 236 $ & $ \sigma  e_{1}\wedge  e_{2} +   e_{1}\wedge  e_{3} -\tau  e_{2}\wedge  e_{3}  $ & 12321321323123 \\
$ 4 $ &  $ 237 $ & $ 2  e_{1}\wedge  e_{3}  $ & 21321321321323 \\
$ 4 $ &  $ 238 $ & $ -2  e_{1}\wedge  e_{3}  $ & 21321321323123 \\
  \hline 

  \hline
 $ 6 $ &  $ 28 $ & $ 1  + \tau  e_{1}\wedge  e_{2} -\sigma  e_{2}\wedge  e_{3}  $ & 1323 \\ 
$ 6 $ &  $ 32 $ & $ 1  + \tau  e_{1}\wedge  e_{2} + \sigma  e_{2}\wedge  e_{3}  $ & 2313 \\ 
$ 6 $ &  $ 37 $ & $ 1  -\tau  e_{1}\wedge  e_{2} -\sigma  e_{2}\wedge  e_{3}  $ & 3132 \\ 
$ 6 $ &  $ 39 $ & $ 1  -\tau  e_{1}\wedge  e_{2} + \sigma  e_{2}\wedge  e_{3}  $ & 3231 \\ 
$ 6 $ &  $ 81 $ & $ 1  -e_{1}\wedge  e_{2} -e_{1}\wedge  e_{3} + e_{2}\wedge  e_{3}  $ & 213232 \\ 
$ 6 $ &  $ 84 $ & $ 1  + e_{1}\wedge  e_{2} + e_{1}\wedge  e_{3} -e_{2}\wedge  e_{3}  $ & 232312 \\ 
$ 6 $ &  $ 123 $ & $ 1  + e_{1}\wedge  e_{2} -e_{1}\wedge  e_{3} -e_{2}\wedge  e_{3}  $ & 12132321 \\ 
$ 6 $ &  $ 126 $ & $ 1  + e_{1}\wedge  e_{2} -e_{1}\wedge  e_{3} + e_{2}\wedge  e_{3}  $ & 12321323 \\ 
$ 6 $ &  $ 127 $ & $ 1  -e_{1}\wedge  e_{2} + e_{1}\wedge  e_{3} + e_{2}\wedge  e_{3}  $ & 12323121 \\ 
$ 6 $ &  $ 132 $ & $ 1  + e_{1}\wedge  e_{2} + e_{1}\wedge  e_{3} + e_{2}\wedge  e_{3}  $ & 13231232 \\ 
$ 6 $ &  $ 136 $ & $ 1  + \sigma  e_{1}\wedge  e_{2} -\tau  e_{1}\wedge  e_{3}  $ & 21323123 \\ 
$ 6 $ &  $ 138 $ & $ 1  -e_{1}\wedge  e_{2} -e_{1}\wedge  e_{3} -e_{2}\wedge  e_{3}  $ & 23213231 \\ 
$ 6 $ &  $ 142 $ & $ 1  -\sigma  e_{1}\wedge  e_{2} + \tau  e_{1}\wedge  e_{3}  $ & 32132312 \\ 
$ 6 $ &  $ 144 $ & $ 1  -e_{1}\wedge  e_{2} + e_{1}\wedge  e_{3} -e_{2}\wedge  e_{3}  $ & 32312321 \\ 
$ 6 $ &  $ 169 $ & $ 1  -\sigma  e_{1}\wedge  e_{2} -\tau  e_{1}\wedge  e_{3}  $ & 1213213213 \\ 
$ 6 $ &  $ 171 $ & $ 1  + \sigma  e_{1}\wedge  e_{3} + \tau  e_{2}\wedge  e_{3}  $ & 1213213232 \\ 
$ 6 $ &  $ 173 $ & $ 1  -\sigma  e_{1}\wedge  e_{3} + \tau  e_{2}\wedge  e_{3}  $ & 1232132132 \\ 
$ 6 $ &  $ 179 $ & $ 1  + \sigma  e_{1}\wedge  e_{2} + \tau  e_{1}\wedge  e_{3}  $ & 1321323121 \\ 
$ 6 $ &  $ 183 $ & $ 1  + \sigma  e_{1}\wedge  e_{3} -\tau  e_{2}\wedge  e_{3}  $ & 2132132321 \\ 
$ 6 $ &  $ 185 $ & $ 1  -\sigma  e_{1}\wedge  e_{3} -\tau  e_{2}\wedge  e_{3}  $ & 2321321321 \\

  \hline 
\end{tabular} \caption{The first set of the conjugacy classes of $\Spin(H_3)$, the ones with orders of `crystallographic type'  $1, 2, 3, 4, 6$. These conjugacy classes  all contain their own reverses i.e.  inverses. The only normal subgroup consists of the first two conjugacy classes and is $\pm 1$, and the order $4$ conjugacy class consists of pure bivectors.  }     \label{tab:CC_spin1}
\end{table}

\begin{table}
\begin{tabular}{|c|c|c|c|c|c|c}
   \hline
   \textbf{Order}  & \textbf{Number} & \textbf{Element} $\times 2$ & \textbf{Generator path}  \\
  \hline 
    \hline 
 $ 5 $ &  $ 8 $ & $ -\tau  +   e_{1}\wedge  e_{2} -\sigma  e_{1}\wedge  e_{3}  $ & 23 \\
$ 5 $ &  $ 10 $ & $ -\tau  -  e_{1}\wedge  e_{2} + \sigma  e_{1}\wedge  e_{3}  $ & 32 \\
$ 5 $ &  $ 23 $ & $ -\tau  -  e_{1}\wedge  e_{2} -\sigma  e_{1}\wedge  e_{3}  $ & 1231 \\
$ 5 $ &  $ 27 $ & $ -\tau  +   e_{1}\wedge  e_{2} + \sigma  e_{1}\wedge  e_{3}  $ & 1321 \\
$ 5 $ &  $ 65 $ & $ -\tau  -\sigma  e_{1}\wedge  e_{2} -  e_{2}\wedge  e_{3}  $ & 121312 \\
$ 5 $ &  $ 67 $ & $ -\tau  +   e_{1}\wedge  e_{3} + \sigma  e_{2}\wedge  e_{3}  $ & 121323 \\
$ 5 $ &  $ 72 $ & $ -\tau  + \sigma  e_{1}\wedge  e_{2} -  e_{2}\wedge  e_{3}  $ & 123232 \\
$ 5 $ &  $ 75 $ & $ -\tau  -  e_{1}\wedge  e_{3} + \sigma  e_{2}\wedge  e_{3}  $ & 132312 \\
$ 5 $ &  $ 77 $ & $ -\tau  + \sigma  e_{1}\wedge  e_{2} +   e_{2}\wedge  e_{3}  $ & 213121 \\
$ 5 $ &  $ 80 $ & $ -\tau  +   e_{1}\wedge  e_{3} -\sigma  e_{2}\wedge  e_{3}  $ & 213231 \\
$ 5 $ &  $ 85 $ & $ -\tau  -\sigma  e_{1}\wedge  e_{2} +   e_{2}\wedge  e_{3}  $ & 232321 \\
$ 5 $ &  $ 90 $ & $ -\tau  -  e_{1}\wedge  e_{3} -\sigma  e_{2}\wedge  e_{3}  $ & 323121 \\

   \hline 
 $ 5 $ &  $ 34 $ & $ -\sigma  -\tau  e_{1}\wedge  e_{2} -  e_{1}\wedge  e_{3}  $ & 2323 \\
$ 5 $ &  $ 40 $ & $ -\sigma  + \tau  e_{1}\wedge  e_{2} +   e_{1}\wedge  e_{3}  $ & 3232 \\
$ 5 $ &  $ 71 $ & $ -\sigma  + \tau  e_{1}\wedge  e_{2} -  e_{1}\wedge  e_{3}  $ & 123231 \\
$ 5 $ &  $ 76 $ & $ -\sigma  -\tau  e_{1}\wedge  e_{2} +   e_{1}\wedge  e_{3}  $ & 132321 \\
$ 5 $ &  $ 122 $ & $ -\sigma  -  e_{1}\wedge  e_{2} + \tau  e_{2}\wedge  e_{3}  $ & 12132312 \\
$ 5 $ &  $ 135 $ & $ -\sigma  +   e_{1}\wedge  e_{2} -\tau  e_{2}\wedge  e_{3}  $ & 21323121 \\
$ 5 $ &  $ 178 $ & $ -\sigma  +   e_{1}\wedge  e_{2} + \tau  e_{2}\wedge  e_{3}  $ & 1321321323 \\
$ 5 $ &  $ 190 $ & $ -\sigma  -  e_{1}\wedge  e_{2} -\tau  e_{2}\wedge  e_{3}  $ & 3213213231 \\
$ 5 $ &  $ 212 $ & $ -\sigma  -\tau  e_{1}\wedge  e_{3} +   e_{2}\wedge  e_{3}  $ & 121321323123 \\
$ 5 $ &  $ 215 $ & $ -\sigma  + \tau  e_{1}\wedge  e_{3} +   e_{2}\wedge  e_{3}  $ & 132132132132 \\
$ 5 $ &  $ 217 $ & $ -\sigma  -\tau  e_{1}\wedge  e_{3} -  e_{2}\wedge  e_{3}  $ & 213213213213 \\
$ 5 $ &  $ 220 $ & $ -\sigma  + \tau  e_{1}\wedge  e_{3} -  e_{2}\wedge  e_{3}  $ & 232132132312 \\
   \hline 
 $ 10 $ &  $ 22 $ & $ \tau  +   e_{1}\wedge  e_{2} + \sigma  e_{1}\wedge  e_{3}  $ & 1213 \\
$ 10 $ &  $ 24 $ & $ \tau  + \sigma  e_{1}\wedge  e_{2} +   e_{2}\wedge  e_{3}  $ & 1232 \\
$ 10 $ &  $ 25 $ & $ \tau  +   e_{1}\wedge  e_{2} -\sigma  e_{1}\wedge  e_{3}  $ & 1312 \\
$ 10 $ &  $ 29 $ & $ \tau  -  e_{1}\wedge  e_{2} + \sigma  e_{1}\wedge  e_{3}  $ & 2131 \\
$ 10 $ &  $ 33 $ & $ \tau  -\sigma  e_{1}\wedge  e_{2} -  e_{2}\wedge  e_{3}  $ & 2321 \\
$ 10 $ &  $ 35 $ & $ \tau  -  e_{1}\wedge  e_{2} -\sigma  e_{1}\wedge  e_{3}  $ & 3121 \\
$ 10 $ &  $ 69 $ & $ \tau  -  e_{1}\wedge  e_{3} -\sigma  e_{2}\wedge  e_{3}  $ & 123123 \\
$ 10 $ &  $ 74 $ & $ \tau  +   e_{1}\wedge  e_{3} -\sigma  e_{2}\wedge  e_{3}  $ & 132132 \\
$ 10 $ &  $ 79 $ & $ \tau  -  e_{1}\wedge  e_{3} + \sigma  e_{2}\wedge  e_{3}  $ & 213213 \\
$ 10 $ &  $ 87 $ & $ \tau  -\sigma  e_{1}\wedge  e_{2} +   e_{2}\wedge  e_{3}  $ & 312323 \\
$ 10 $ &  $ 88 $ & $ \tau  +   e_{1}\wedge  e_{3} + \sigma  e_{2}\wedge  e_{3}  $ & 321321 \\
$ 10 $ &  $ 92 $ & $ \tau  + \sigma  e_{1}\wedge  e_{2} -  e_{2}\wedge  e_{3}  $ & 323213 \\

 \hline 
 $ 10 $ &  $ 70 $ & $ \sigma  -\tau  e_{1}\wedge  e_{2} +   e_{1}\wedge  e_{3}  $ & 123213 \\
$ 10 $ &  $ 73 $ & $ \sigma  -\tau  e_{1}\wedge  e_{2} -  e_{1}\wedge  e_{3}  $ & 131232 \\
$ 10 $ &  $ 78 $ & $ \sigma  + \tau  e_{1}\wedge  e_{2} +   e_{1}\wedge  e_{3}  $ & 213123 \\
$ 10 $ &  $ 86 $ & $ \sigma  + \tau  e_{1}\wedge  e_{2} -  e_{1}\wedge  e_{3}  $ & 312321 \\
$ 10 $ &  $ 121 $ & $ \sigma  +   e_{1}\wedge  e_{2} -\tau  e_{2}\wedge  e_{3}  $ & 12132132 \\
$ 10 $ &  $ 133 $ & $ \sigma  -  e_{1}\wedge  e_{2} + \tau  e_{2}\wedge  e_{3}  $ & 21321321 \\
$ 10 $ &  $ 180 $ & $ \sigma  -  e_{1}\wedge  e_{2} -\tau  e_{2}\wedge  e_{3}  $ & 1321323123 \\
$ 10 $ &  $ 189 $ & $ \sigma  +   e_{1}\wedge  e_{2} + \tau  e_{2}\wedge  e_{3}  $ & 3213213213 \\
$ 10 $ &  $ 210 $ & $ \sigma  + \tau  e_{1}\wedge  e_{3} -  e_{2}\wedge  e_{3}  $ & 121321321323 \\
$ 10 $ &  $ 216 $ & $ \sigma  -\tau  e_{1}\wedge  e_{3} -  e_{2}\wedge  e_{3}  $ & 132132132312 \\
$ 10 $ &  $ 218 $ & $ \sigma  + \tau  e_{1}\wedge  e_{3} +   e_{2}\wedge  e_{3}  $ & 213213213231 \\
$ 10 $ &  $ 219 $ & $ \sigma  -\tau  e_{1}\wedge  e_{3} +   e_{2}\wedge  e_{3}  $ & 232132132132 \\

  \hline 
\end{tabular}\caption{The second set of  conjugacy classes of $\Spin(H_3)$, the ones with orders of `non-crystallographic type' i.e. the ones related to 5-fold symmetry. These conjugacy classes also all contain their inverses. The first column denotes the order of the elements in each conjugacy class. The second column is the position in the order in which our algorithm generates this element, for convenience (c.f. the supplementary material). The final column denotes the order in which the generators with the corresponding labels are applied to generate this group element, i.e. is effectively the `word in the generators' that yields this element. To avoid confusion, pairs of such words doubly cover the rotations of the icosahedral group $A_5$, which are also often considered in terms of words in the $A_5$ generators. The ones meant here are the root vectors multiplied by using the geometric product.  }     \label{tab:CC_spin2}
\end{table}
\normalsize

The spin group $\Spin(3)=2I$ is given in Tables \ref{tab:CC_spin1} to \ref{tab:CC_spin2}. Its nine conjugacy classes lead to irreducible representations of dimensions $1,3,3,4,5$ which are shared by $A_5$, as well as the spinorial ones of dimensions $2,2,4,6$. (Interesting connections with the binary polyhedral groups and the McKay correspondence \cite{Mckay1980graphs} are explored elsewhere \cite{Dechant2018Trinity}). Since $\Pin(H_3) = \Spin(H_3) \times  \mathbb{Z}_2$, the remaining Tables \ref{tab:CC_pin1} to \ref{tab:CC_pin2} list the remaining 9 conjugacy classes achieved by multiplying those of $\Spin(H_3)$ with the inversion $e_1e_2e_3$. 

Tables \ref{tab:CC_spin1} to \ref{tab:CC_spin2}: The only normal subgroup of $2I$ consists of the first two conjugacy classes i.e. $\pm 1$. We note that the conjugacy class of order $4$ consists of thirty pure bivectors, and that they give rise to the  2-fold rotations of the icosahedron around its $30$ edges. The $20$ 3-fold rotations around the $20$ triangular faces are split into two conjugacy classes of orders $3$ and $6$, which are related by multiplication by $-1$, and each contain their own reverse/inverse.  Similarly, the two sets of twelve 5-fold rotations are doubly covered by 4 conjugacy classes which are related by multiplication by $-1$. The two classes describe rotations by $\pm 2\pi/5$ and $\pm 4\pi/5$, respectively, around the 5-fold axes of symmetry, the icosahedral vertices. 

Tables \ref{tab:CC_pin1} to \ref{tab:CC_pin2}: The first two conjugacy classes are the inversion and its negative. The conjugacy class consisting of pure vectors of course corresponds to the $30$ roots of $H_3$ which generate the reflections, and which are of course related to the $30$ 2-fold rotations since the inversion is contained in the group (so one can dualise a (root) vector to a pure bivector). The conjugacy classes of order $12$ are the two inversion-related versions of the 3-fold rotations, and are rotoreflections. The four conjugacy classes of order $20$ are both related to the 5-fold rotations, as well as serving as the versor analogues of the Coxeter elements e.g. $w=a_1a_2a_3$. These are in one conjugacy class in the  reflection/Coxeter group framework where their order gives the Coxeter number. But in this Clifford double cover setup these `Coxeter versors' are given in 4 conjugacy classes that are related by reversal and multiplication by $e_1e_2e_3$. 

\begin{rem}
It has been noted that $a_1a_2$ and $a_2a_3$ generate the quaternionic root system multiplicatively e.g. for $H_4$. This is pretty obvious when thought of in terms of the 3D simple roots and the Induction Theorem, as they of course generate $\Spin(H_3)$, which gives rise to the $H_4$ root system. 
\end{rem}

This example illustrates how one can perform practical computations in group theory via versors in this Clifford algebra framework, and in \verb|galgebra| in particular. We will discuss group and representation theoretic aspects in more detail elsewhere.

\tiny

\begin{table}
\begin{tabular}{|c|c|c|c|c|c|c}
   \hline
   \textbf{Order}  & \textbf{Number} & \textbf{Element} $\times 2$ & \textbf{Generator path}  \\
  \hline 

 $ 4 $ &  $ 240 $ & $ 2  e_{1}\wedge  e_{2}\wedge  e_{3}  $ & 121321321323123 \\

  \hline 
 $ 4 $ &  $ 239 $ & $ -2  e_{1}\wedge  e_{2}\wedge  e_{3}  $ & 121321321321323 \\

   \hline 
 $ 2 $ &  $ 1 $ & $ 2  e_{2}  $ & 1 \\
$ 2 $ &  $ 2 $ & $ -\tau  e_{1} -  e_{2} + \sigma  e_{3}  $ & 2 \\
$ 2 $ &  $ 3 $ & $ 2  e_{1}  $ & 3 \\
$ 2 $ &  $ 11 $ & $ \tau  e_{1} -  e_{2} -\sigma  e_{3}  $ & 121 \\
$ 2 $ &  $ 13 $ & $ -2  e_{1}  $ & 131 \\
$ 2 $ &  $ 17 $ & $ -\sigma  e_{1} + \tau  e_{2} +   e_{3}  $ & 232 \\
$ 2 $ &  $ 19 $ & $ -2  e_{2}  $ & 313 \\
$ 2 $ &  $ 21 $ & $ -\tau  e_{1} +   e_{2} -\sigma  e_{3}  $ & 323 \\
$ 2 $ &  $ 44 $ & $ -\tau  e_{1} +   e_{2} + \sigma  e_{3}  $ & 12313 \\
$ 2 $ &  $ 45 $ & $ \sigma  e_{1} + \tau  e_{2} -  e_{3}  $ & 12321 \\
$ 2 $ &  $ 48 $ & $ \tau  e_{1} -  e_{2} + \sigma  e_{3}  $ & 13123 \\
$ 2 $ &  $ 49 $ & $ \tau  e_{1} +   e_{2} -\sigma  e_{3}  $ & 13132 \\
$ 2 $ &  $ 50 $ & $ -\tau  e_{1} -  e_{2} -\sigma  e_{3}  $ & 13213 \\
$ 2 $ &  $ 51 $ & $ \tau  e_{1} +   e_{2} + \sigma  e_{3}  $ & 13231 \\
$ 2 $ &  $ 53 $ & $ \sigma  e_{1} -\tau  e_{2} -  e_{3}  $ & 21312 \\
$ 2 $ &  $ 60 $ & $ -\sigma  e_{1} -\tau  e_{2} -  e_{3}  $ & 23232 \\
$ 2 $ &  $ 93 $ & $ -\sigma  e_{1} -\tau  e_{2} +   e_{3}  $ & 1213121 \\
$ 2 $ &  $ 101 $ & $ \sigma  e_{1} -\tau  e_{2} +   e_{3}  $ & 1232321 \\
$ 2 $ &  $ 103 $ & $ \sigma  e_{1} + \tau  e_{2} +   e_{3}  $ & 1312323 \\
$ 2 $ &  $ 108 $ & $ -\sigma  e_{1} + \tau  e_{2} -  e_{3}  $ & 1323213 \\
$ 2 $ &  $ 109 $ & $ -  e_{1} + \sigma  e_{2} -\tau  e_{3}  $ & 2132132 \\
$ 2 $ &  $ 110 $ & $   e_{1} -\sigma  e_{2} + \tau  e_{3}  $ & 2132312 \\
$ 2 $ &  $ 145 $ & $   e_{1} + \sigma  e_{2} + \tau  e_{3}  $ & 121321321 \\
$ 2 $ &  $ 147 $ & $ -  e_{1} -\sigma  e_{2} -\tau  e_{3}  $ & 121323121 \\
$ 2 $ &  $ 166 $ & $ -  e_{1} -\sigma  e_{2} + \tau  e_{3}  $ & 321321323 \\
$ 2 $ &  $ 168 $ & $   e_{1} + \sigma  e_{2} -\tau  e_{3}  $ & 321323123 \\
$ 2 $ &  $ 199 $ & $ -  e_{1} + \sigma  e_{2} + \tau  e_{3}  $ & 13213213213 \\
$ 2 $ &  $ 200 $ & $   e_{1} -\sigma  e_{2} -\tau  e_{3}  $ & 13213213231 \\
$ 2 $ &  $ 229 $ & $ -2  e_{3}  $ & 2132132132132 \\
$ 2 $ &  $ 230 $ & $ 2  e_{3}  $ & 2132132132312 \\

\hline 
 $ 12 $ &  $ 52 $ & $ -\tau  e_{1} -\sigma  e_{2} -  e_{1}\wedge  e_{2}\wedge  e_{3}  $ & 13232 \\
$ 12 $ &  $ 55 $ & $ -\tau  e_{2} + \sigma  e_{3} -  e_{1}\wedge  e_{2}\wedge  e_{3}  $ & 21323 \\
$ 12 $ &  $ 58 $ & $ \tau  e_{1} + \sigma  e_{2} -  e_{1}\wedge  e_{2}\wedge  e_{3}  $ & 23213 \\
$ 12 $ &  $ 62 $ & $ \tau  e_{2} -\sigma  e_{3} -  e_{1}\wedge  e_{2}\wedge  e_{3}  $ & 32132 \\
$ 12 $ &  $ 64 $ & $ \tau  e_{1} -\sigma  e_{2} -  e_{1}\wedge  e_{2}\wedge  e_{3}  $ & 32321 \\
$ 12 $ &  $ 94 $ & $ -\tau  e_{1} + \sigma  e_{2} -  e_{1}\wedge  e_{2}\wedge  e_{3}  $ & 1213123 \\
$ 12 $ &  $ 96 $ & $ -\tau  e_{2} -\sigma  e_{3} -  e_{1}\wedge  e_{2}\wedge  e_{3}  $ & 1213231 \\
$ 12 $ &  $ 98 $ & $ -  e_{1} -  e_{2} -  e_{3} -  e_{1}\wedge  e_{2}\wedge  e_{3}  $ & 1231232 \\
$ 12 $ &  $ 100 $ & $ -  e_{1} +   e_{2} -  e_{3} -  e_{1}\wedge  e_{2}\wedge  e_{3}  $ & 1232312 \\
$ 12 $ &  $ 104 $ & $ \tau  e_{2} + \sigma  e_{3} -  e_{1}\wedge  e_{2}\wedge  e_{3}  $ & 1321321 \\
$ 12 $ &  $ 112 $ & $   e_{1} -  e_{2} +   e_{3} -  e_{1}\wedge  e_{2}\wedge  e_{3}  $ & 2312321 \\
$ 12 $ &  $ 115 $ & $   e_{1} +   e_{2} +   e_{3} -  e_{1}\wedge  e_{2}\wedge  e_{3}  $ & 2323121 \\
$ 12 $ &  $ 116 $ & $   e_{1} -  e_{2} -  e_{3} -  e_{1}\wedge  e_{2}\wedge  e_{3}  $ & 2323123 \\
$ 12 $ &  $ 120 $ & $ -  e_{1} +   e_{2} +   e_{3} -  e_{1}\wedge  e_{2}\wedge  e_{3}  $ & 3231232 \\
$ 12 $ &  $ 149 $ & $ -  e_{1} -  e_{2} +   e_{3} -  e_{1}\wedge  e_{2}\wedge  e_{3}  $ & 123213213 \\
$ 12 $ &  $ 156 $ & $   e_{1} +   e_{2} -  e_{3} -  e_{1}\wedge  e_{2}\wedge  e_{3}  $ & 132312321 \\
$ 12 $ &  $ 191 $ & $ \sigma  e_{1} -\tau  e_{3} -  e_{1}\wedge  e_{2}\wedge  e_{3}  $ & 12132132132 \\
$ 12 $ &  $ 201 $ & $ -\sigma  e_{1} + \tau  e_{3} -  e_{1}\wedge  e_{2}\wedge  e_{3}  $ & 21321321321 \\
$ 12 $ &  $ 228 $ & $ \sigma  e_{1} + \tau  e_{3} -  e_{1}\wedge  e_{2}\wedge  e_{3}  $ & 1321321323123 \\
$ 12 $ &  $ 232 $ & $ -\sigma  e_{1} -\tau  e_{3} -  e_{1}\wedge  e_{2}\wedge  e_{3}  $ & 2321321323123 \\
\hline 
 $ 12 $ &  $ 46 $ & $ \tau  e_{1} -\sigma  e_{2} +   e_{1}\wedge  e_{2}\wedge  e_{3}  $ & 12323 \\
$ 12 $ &  $ 57 $ & $ \tau  e_{2} -\sigma  e_{3} +   e_{1}\wedge  e_{2}\wedge  e_{3}  $ & 23123 \\
$ 12 $ &  $ 59 $ & $ -\tau  e_{1} -\sigma  e_{2} +   e_{1}\wedge  e_{2}\wedge  e_{3}  $ & 23231 \\
$ 12 $ &  $ 61 $ & $ \tau  e_{1} + \sigma  e_{2} +   e_{1}\wedge  e_{2}\wedge  e_{3}  $ & 31232 \\
$ 12 $ &  $ 63 $ & $ -\tau  e_{2} + \sigma  e_{3} +   e_{1}\wedge  e_{2}\wedge  e_{3}  $ & 32312 \\
$ 12 $ &  $ 95 $ & $ \tau  e_{2} + \sigma  e_{3} +   e_{1}\wedge  e_{2}\wedge  e_{3}  $ & 1213213 \\
$ 12 $ &  $ 97 $ & $   e_{1} +   e_{2} +   e_{3} +   e_{1}\wedge  e_{2}\wedge  e_{3}  $ & 1213232 \\
$ 12 $ &  $ 99 $ & $   e_{1} -  e_{2} +   e_{3} +   e_{1}\wedge  e_{2}\wedge  e_{3}  $ & 1232132 \\
$ 12 $ &  $ 102 $ & $ -\tau  e_{1} + \sigma  e_{2} +   e_{1}\wedge  e_{2}\wedge  e_{3}  $ & 1312321 \\
$ 12 $ &  $ 106 $ & $ -\tau  e_{2} -\sigma  e_{3} +   e_{1}\wedge  e_{2}\wedge  e_{3}  $ & 1323121 \\
$ 12 $ &  $ 111 $ & $ -  e_{1} +   e_{2} -  e_{3} +   e_{1}\wedge  e_{2}\wedge  e_{3}  $ & 2132321 \\
$ 12 $ &  $ 113 $ & $ -  e_{1} -  e_{2} -  e_{3} +   e_{1}\wedge  e_{2}\wedge  e_{3}  $ & 2321321 \\
$ 12 $ &  $ 114 $ & $ -  e_{1} +   e_{2} +   e_{3} +   e_{1}\wedge  e_{2}\wedge  e_{3}  $ & 2321323 \\
$ 12 $ &  $ 119 $ & $   e_{1} -  e_{2} -  e_{3} +   e_{1}\wedge  e_{2}\wedge  e_{3}  $ & 3213232 \\
$ 12 $ &  $ 150 $ & $   e_{1} +   e_{2} -  e_{3} +   e_{1}\wedge  e_{2}\wedge  e_{3}  $ & 123213231 \\
$ 12 $ &  $ 155 $ & $ -  e_{1} -  e_{2} +   e_{3} +   e_{1}\wedge  e_{2}\wedge  e_{3}  $ & 132132321 \\
$ 12 $ &  $ 192 $ & $ -\sigma  e_{1} + \tau  e_{3} +   e_{1}\wedge  e_{2}\wedge  e_{3}  $ & 12132132312 \\
$ 12 $ &  $ 203 $ & $ \sigma  e_{1} -\tau  e_{3} +   e_{1}\wedge  e_{2}\wedge  e_{3}  $ & 21321323121 \\
$ 12 $ &  $ 227 $ & $ -\sigma  e_{1} -\tau  e_{3} +   e_{1}\wedge  e_{2}\wedge  e_{3}  $ & 1321321321323 \\
$ 12 $ &  $ 231 $ & $ \sigma  e_{1} + \tau  e_{3} +   e_{1}\wedge  e_{2}\wedge  e_{3}  $ & 2321321321323 \\
 \hline

  \hline 
\end{tabular} \caption{The first set of the remaining conjugacy classes of $\Pin(H_3)$. Note that the first two conjugacy classes show that the inversion $e_1e_2e_3$ is contained in the group. The class of order $2$ consists of pure vectors and are thus their own reverse i.e. inverse. The other pair are each other's reverses/inverses. }     \label{tab:CC_pin1}
\end{table}

\begin{table}
\begin{tabular}{|c|c|c|c|c|c|c}
   \hline
   \textbf{Order}  & \textbf{Number} & \textbf{Element} $\times 2$ & \textbf{Generator path}  \\
  \hline 
    \hline 

 $ 20 $ &  $ 14 $ & $   e_{1} -\tau  e_{2} -\sigma  e_{1}\wedge  e_{2}\wedge  e_{3}  $ & 132 \\
$ 20 $ &  $ 15 $ & $ -  e_{1} + \tau  e_{2} -\sigma  e_{1}\wedge  e_{2}\wedge  e_{3}  $ & 213 \\
$ 20 $ &  $ 20 $ & $ -  e_{1} -\tau  e_{2} -\sigma  e_{1}\wedge  e_{2}\wedge  e_{3}  $ & 321 \\
$ 20 $ &  $ 41 $ & $   e_{1} + \tau  e_{2} -\sigma  e_{1}\wedge  e_{2}\wedge  e_{3}  $ & 12131 \\
$ 20 $ &  $ 43 $ & $ \tau  e_{1} +   e_{3} -\sigma  e_{1}\wedge  e_{2}\wedge  e_{3}  $ & 12312 \\
$ 20 $ &  $ 56 $ & $ -\tau  e_{1} -  e_{3} -\sigma  e_{1}\wedge  e_{2}\wedge  e_{3}  $ & 23121 \\
$ 20 $ &  $ 105 $ & $ \tau  e_{1} -  e_{3} -\sigma  e_{1}\wedge  e_{2}\wedge  e_{3}  $ & 1321323 \\
$ 20 $ &  $ 118 $ & $ -\tau  e_{1} +   e_{3} -\sigma  e_{1}\wedge  e_{2}\wedge  e_{3}  $ & 3213231 \\
$ 20 $ &  $ 159 $ & $   e_{2} + \tau  e_{3} -\sigma  e_{1}\wedge  e_{2}\wedge  e_{3}  $ & 213213232 \\
$ 20 $ &  $ 162 $ & $ -  e_{2} -\tau  e_{3} -\sigma  e_{1}\wedge  e_{2}\wedge  e_{3}  $ & 232132312 \\
$ 20 $ &  $ 193 $ & $   e_{2} -\tau  e_{3} -\sigma  e_{1}\wedge  e_{2}\wedge  e_{3}  $ & 12132132321 \\
$ 20 $ &  $ 197 $ & $ -  e_{2} + \tau  e_{3} -\sigma  e_{1}\wedge  e_{2}\wedge  e_{3}  $ & 12321323121 \\
  \hline 
 $ 20 $ &  $ 12 $ & $ -  e_{1} -\tau  e_{2} + \sigma  e_{1}\wedge  e_{2}\wedge  e_{3}  $ & 123 \\
$ 20 $ &  $ 16 $ & $   e_{1} -\tau  e_{2} + \sigma  e_{1}\wedge  e_{2}\wedge  e_{3}  $ & 231 \\
$ 20 $ &  $ 18 $ & $ -  e_{1} + \tau  e_{2} + \sigma  e_{1}\wedge  e_{2}\wedge  e_{3}  $ & 312 \\
$ 20 $ &  $ 42 $ & $ -\tau  e_{1} -  e_{3} + \sigma  e_{1}\wedge  e_{2}\wedge  e_{3}  $ & 12132 \\
$ 20 $ &  $ 47 $ & $   e_{1} + \tau  e_{2} + \sigma  e_{1}\wedge  e_{2}\wedge  e_{3}  $ & 13121 \\
$ 20 $ &  $ 54 $ & $ \tau  e_{1} +   e_{3} + \sigma  e_{1}\wedge  e_{2}\wedge  e_{3}  $ & 21321 \\
$ 20 $ &  $ 107 $ & $ -\tau  e_{1} +   e_{3} + \sigma  e_{1}\wedge  e_{2}\wedge  e_{3}  $ & 1323123 \\
$ 20 $ &  $ 117 $ & $ \tau  e_{1} -  e_{3} + \sigma  e_{1}\wedge  e_{2}\wedge  e_{3}  $ & 3213213 \\
$ 20 $ &  $ 160 $ & $ -  e_{2} -\tau  e_{3} + \sigma  e_{1}\wedge  e_{2}\wedge  e_{3}  $ & 213231232 \\
$ 20 $ &  $ 161 $ & $   e_{2} + \tau  e_{3} + \sigma  e_{1}\wedge  e_{2}\wedge  e_{3}  $ & 232132132 \\
$ 20 $ &  $ 194 $ & $ -  e_{2} + \tau  e_{3} + \sigma  e_{1}\wedge  e_{2}\wedge  e_{3}  $ & 12132312321 \\
$ 20 $ &  $ 195 $ & $   e_{2} -\tau  e_{3} + \sigma  e_{1}\wedge  e_{2}\wedge  e_{3}  $ & 12321321321 \\
  \hline 
 $ 20 $ &  $ 146 $ & $ \sigma  e_{1} -  e_{2} -\tau  e_{1}\wedge  e_{2}\wedge  e_{3}  $ & 121321323 \\
$ 20 $ &  $ 151 $ & $ -  e_{1} + \sigma  e_{3} -\tau  e_{1}\wedge  e_{2}\wedge  e_{3}  $ & 123213232 \\
$ 20 $ &  $ 154 $ & $ \sigma  e_{1} +   e_{2} -\tau  e_{1}\wedge  e_{2}\wedge  e_{3}  $ & 132132312 \\
$ 20 $ &  $ 158 $ & $ -\sigma  e_{1} -  e_{2} -\tau  e_{1}\wedge  e_{2}\wedge  e_{3}  $ & 213213231 \\
$ 20 $ &  $ 163 $ & $   e_{1} -\sigma  e_{3} -\tau  e_{1}\wedge  e_{2}\wedge  e_{3}  $ & 232132321 \\
$ 20 $ &  $ 167 $ & $ -\sigma  e_{1} +   e_{2} -\tau  e_{1}\wedge  e_{2}\wedge  e_{3}  $ & 321323121 \\
$ 20 $ &  $ 198 $ & $ -  e_{1} -\sigma  e_{3} -\tau  e_{1}\wedge  e_{2}\wedge  e_{3}  $ & 12321323123 \\
$ 20 $ &  $ 202 $ & $ \sigma  e_{2} -  e_{3} -\tau  e_{1}\wedge  e_{2}\wedge  e_{3}  $ & 21321321323 \\
$ 20 $ &  $ 205 $ & $   e_{1} + \sigma  e_{3} -\tau  e_{1}\wedge  e_{2}\wedge  e_{3}  $ & 23213213213 \\
$ 20 $ &  $ 207 $ & $ -\sigma  e_{2} +   e_{3} -\tau  e_{1}\wedge  e_{2}\wedge  e_{3}  $ & 32132132132 \\
$ 20 $ &  $ 224 $ & $ \sigma  e_{2} +   e_{3} -\tau  e_{1}\wedge  e_{2}\wedge  e_{3}  $ & 1213213213231 \\
$ 20 $ &  $ 226 $ & $ -\sigma  e_{2} -  e_{3} -\tau  e_{1}\wedge  e_{2}\wedge  e_{3}  $ & 1232132132312 \\
  \hline 
 $ 20 $ &  $ 148 $ & $ -\sigma  e_{1} +   e_{2} + \tau  e_{1}\wedge  e_{2}\wedge  e_{3}  $ & 121323123 \\
$ 20 $ &  $ 152 $ & $   e_{1} -\sigma  e_{3} + \tau  e_{1}\wedge  e_{2}\wedge  e_{3}  $ & 123231232 \\
$ 20 $ &  $ 153 $ & $ -\sigma  e_{1} -  e_{2} + \tau  e_{1}\wedge  e_{2}\wedge  e_{3}  $ & 132132132 \\
$ 20 $ &  $ 157 $ & $ \sigma  e_{1} +   e_{2} + \tau  e_{1}\wedge  e_{2}\wedge  e_{3}  $ & 213213213 \\
$ 20 $ &  $ 164 $ & $ -  e_{1} + \sigma  e_{3} + \tau  e_{1}\wedge  e_{2}\wedge  e_{3}  $ & 232312321 \\
$ 20 $ &  $ 165 $ & $ \sigma  e_{1} -  e_{2} + \tau  e_{1}\wedge  e_{2}\wedge  e_{3}  $ & 321321321 \\
$ 20 $ &  $ 196 $ & $   e_{1} + \sigma  e_{3} + \tau  e_{1}\wedge  e_{2}\wedge  e_{3}  $ & 12321321323 \\
$ 20 $ &  $ 204 $ & $ -\sigma  e_{2} +   e_{3} + \tau  e_{1}\wedge  e_{2}\wedge  e_{3}  $ & 21321323123 \\
$ 20 $ &  $ 206 $ & $ -  e_{1} -\sigma  e_{3} + \tau  e_{1}\wedge  e_{2}\wedge  e_{3}  $ & 23213213231 \\
$ 20 $ &  $ 208 $ & $ \sigma  e_{2} -  e_{3} + \tau  e_{1}\wedge  e_{2}\wedge  e_{3}  $ & 32132132312 \\
$ 20 $ &  $ 223 $ & $ -\sigma  e_{2} -  e_{3} + \tau  e_{1}\wedge  e_{2}\wedge  e_{3}  $ & 1213213213213 \\
$ 20 $ &  $ 225 $ & $ \sigma  e_{2} +   e_{3} + \tau  e_{1}\wedge  e_{2}\wedge  e_{3}  $ & 1232132132132 \\

  \hline 
\end{tabular}\caption{The four conjugacy classes of order $20$ are  related to the 5-fold rotations in $\Spin(H_3)$ via multiplication by $e_1e_2e_3$. They are also the spinor version of the  Coxeter elements e.g. $w=a_1a_2a_3$. They are related by reversal and multiplication by $e_1e_2e_3$.  }     \label{tab:CC_pin2}
\end{table}

\normalsize

\section{Subgroups}\label{sec_group}
The Induction Theorem \ref{thrm_ind} from Section \ref{sec_background} showed that every 3D root system determines a 4D root system. This proceeded essentially via using the 3D roots to construct a group of spinors (via multiplication with the geometric product), which satisfies the properties of a root system. Our main example is of course $H_3$ which induces $H_4$ in four dimensions via the  binary icosahedral group of order $120$ as the spinor group intermediary. It is therefore a straightforward corollary of the Induction Theorem that each even subgroup of $2I$ also yields a root system.

\begin{cor}[Subgroups of $2I$]
Each even subgroup $G$ of the binary icosahedral group $2I$ determines a corresponding root system $\Phi$ that is a subset of the $H_4$ root system, the 600-cell. 
\end{cor}

\begin{thm}[Induced subrootsystems of $H_4$]
The binary icosahedral group $2I$ has the following subgroups that determines the corresponding root systems:
\begin{itemize}
\item The normal subgroup $\pm 1$ which gives $A_1$.
\item The quaternion group $Q$ consisting of $\pm 1$, $\pm e_1 e_2$, $\pm e_2 e_3$ and $\pm e_3 e_1$, which gives $A_1\times A_1\times A_1\times A_1$.
\item The binary dihedral groups of orders $6$ and $10$, which yield $A_2$ and $H_2$.
\item The binary tetrahedral group, which yields $D_4$. 
\end{itemize}
\end{thm}
\begin{rem}
Note that although the $A_3$ root system is not contained in $H_3$, $2T$ is a subgroup of $2I$ and therefore $D_4$ is contained in $H_4$. We will revisit these examples in later sections and in the next section investigate this delicate relationship between subgroups and other subrootsystems further. 
\end{rem}

\begin{prop}[Simple roots of induced subrootsystems]
$A_2$ and $H_2$ are generated straightforwardly from the $H_3$ generators $a_1, a_2, a_3$ e.g. via the `spinorial simple roots' $a_1a_1=1$ and $a_1a_2$ for $A_2$ and $a_1a_1=1$ and $a_2a_3$ for $H_2$. 

One possible choice of simple roots for $D_4$ contained in $H_4$ is given by $$(a_1a_1, a_1a_2, a_1a_2a_3a_2a_3a_1a_2a_3, a_3a_2a_1a_3a_2a_1a_3a_2),$$ but it is of course not unique. 
\end{prop}
Explicit versions of these simple roots can be looked up in the earlier tables via the `generator path', which can be used to explicitly verify the correct Cartan matrix and closure of the root system (see e.g. supplementary information).

\section{Subrootsystems}\label{sec_rs}
There is a subtlety at play here since there are subrootsystems of $H_4$ that are neither induced by 3D subrootsystems nor by even spinor subgroups. We have observed in the spinor Induction Theorem that $A_1\oplus I_2(n)$ root systems experience a doubling to $I_2(n)\oplus I_2(n)$ in the induction process. Indeed, such $A_2\oplus A_2$ and $H_2\oplus H_2$ within $H_4$ are induced, but not because they have an orthogonal $A_1$. Instead, it is because the group $H_3$ contains the inversion, which is manifested at the level of the pin group by virtue of containing the pseudoscalar $e_1e_2e_3$ from Table \ref{tab:CC_pin1}. These are all contained in $H_3$, and have the effect of creating a second orthogonal $I_2(n)$ in the even subalgebra. We illustrate the idea and its complexities with some examples. 
\begin{ex}
Take as a first example $A_1^3$ with simple roots $e_1, e_2, e_3$ but think of it as $I_2(2)\oplus A_1$. Now this gets doubled to $I_2(2)\oplus I_2(2)= A_1^4$ via the spinor group $\pm 1$, $\pm e_1 e_2$, $\pm e_2 e_3$ and $\pm e_3 e_1$. \end{ex}
\begin{ex}
Now consider the following twist: take $A_1^2$ with simple roots $e_1, e_2$ but instead of having $e_3$ available as another orthogonal simple root, we just have the inversion $e_1e_2e_3$ available. So we get $\pm 1$ and $\pm e_1e_2$ from multiplying the simple roots. But by operating in the whole pin group we can multiply $e_1$ by the pseudoscalar  $e_1e_2e_3$, which yields $e_2e_3$, which \emph{is} in the spin part. So similarly we get $\pm e_2 e_3$ and $\pm e_3 e_1$, i.e. we get the same spinor group as in the previous example, without actually having the third simple root $e_3$ available. This therefore induces the same $A_1^4$ root system. 
\end{ex}

\begin{prop}[Doubling -- even case]
For even $n$ the root system $I_2(n)$ together with the inversion $e_1e_2e_3$ yields the doubling $I_2(n)\oplus I_2(n)$.  \end{prop}

\begin{proof}
Without loss of generality take $e_1$ as the first simple root. Since $n$ is even, the number of roots is a multiple of $4$ and therefore $e_2$ is also a root. Therefore having $e_1$, $e_2$ and  $e_1e_2e_3$ available is equivalent to having $e_3$ available  as well, which via the Induction Theorem leads to a doubling $I_2(n)\oplus I_2(n)$.
\end{proof}

A convenient choice of simple roots for $I_2(n)$ is  $\alpha_1=e_1$ and $\alpha_2=-\cos{\frac{\pi}{n}}e_1+\sin{\frac{\pi}{n}}e_2$.

\begin{prop}[Doubling -- general case]
The root system $I_2(n)$ together with the inversion $e_1e_2e_3$ yields the doubling $I_2(n)\oplus I_2(n)$. A possible choice of simple roots for $I_2(n)\oplus I_2(n)$ is given by $\alpha_1\alpha_1, \alpha_1\alpha_2, \alpha_1e_1e_2e_3, \alpha_2e_1e_2e_3)$. \end{prop}

\begin{proof}
Direct calculation confirms that these constitute two orthogonal $I_2(n)$ root systems with respect to the spinor inner product, and the simple roots give the correct Cartan matrix
$$A \left(I_2(n)\oplus I_2(n)\right) = \begin{pmatrix}
	    2&-2\cos\frac{\pi}{n}&0&0
	\\ -2\cos\frac{\pi}{n}&2&0&0 
	\\  0&0&2&-2\cos\frac{\pi}{n}
	\\  0&0&-2\cos\frac{\pi}{n}&2
	 \end{pmatrix}.$$
C.f. also  the computational proof in the supplementary information. 
\end{proof}

\begin{rem}
The inclusion of the $A_2\oplus A_2$ or $H_2\oplus H_2$ subrootsystems therefore arises in the more oblique way for odd $n$ root systems as $n=3$ and $n=5$. The Cartan matrices are given by $$A \left(A_2\oplus A_2\right) = \begin{pmatrix}
	    2&-1&0&0
	\\ -1&2&0&0 
	\\  0&0&2&-1
	\\  0&0&-1&2
	 \end{pmatrix}, A \left(H_2\oplus H_2\right) = \begin{pmatrix}
	    2&-\tau&0&0
	\\ -\tau&2&0&0 
	\\  0&0&2&-\tau
	\\  0&0&-\tau&2
	 \end{pmatrix}$$ \end{rem}

\begin{prop}[Simple roots of other subrootsystems]
$A_2\oplus A_2$ and $H_2\oplus H_2$ are generated straightforwardly from the $H_3$ generators $a_1, a_2, a_3$ e.g. via the `spinorial simple roots' $$a_1a_1=1, a_1a_2, a_1e_1e_2e_3, a_2e_1e_2e_3$$ for $A_2$ and $$a_1a_1=1, a_2a_3, a_2e_1e_2e_3, a_3e_1e_2e_3$$for $H_2$. 

One possible choice of simple roots for $A_4$ within $H_4$ is given by $$(a_1a_1, a_1a_2, a_1a_3a_2a_1a_3a_2a_1a_3, a_3a_2a_1a_3a_2a_1a_3a_2a_3a_1a_2a_3),$$ but other choices are of course possible. 
\end{prop}

\begin{rem}
The fact that $H_3$ can't have $A_1\oplus A_2$ or $A_1\oplus H_2$ subrootsystems is clear from the following: the $H_3$ root system is the icosidodecahedron with vertices at the 2-fold axes. It contains decagonal/hexagonal grand circles which are valid   $A_2$ or $H_2$ subrootsystems. However, a root normal to those can't exist because they would be the vertices of the icosahedron (5-fold axes) or dodecahedron (3-fold axes), which is of course inconsistent with the vertices being the 2-fold axes of the icosidodecahedron.  \end{rem}

\begin{rem}
For odd $n$ the 4D root systems induced by $A_1\oplus I_2(n)$ and via  $I_2(n)$ in combination with the  inversion $e_1e_2e_3$ are subtly different (related via $e_1 \leftrightarrow e_2$). However, for even $n$ they of course coincide. \end{rem}

\begin{rem}
The inversion $e_1e_2e_3$ is often contained in Coxeter groups, but is famously not contained in the $A_n$ family for odd $n$. As such $A_3$, the tetrahedral group, is a prime example of where this isn't the case; it is even obvious from the tetrahedron itself that it is not inversion invariant. The existence of the inversion in a group means that one can use this pseudoscalar to dualise root vectors to pure bivectors. In work on quaternions it was often regarded as deeply meaningful that the pure quaternion roots e.g. of $H_4$ are exactly the $H_3$ root system, and analogously for other cases (e.g. \cite{moody1993quasicrystals}). However, a more useful way of viewing this is that this is pretty obvious from the Induction Theorem as long as the inversion is contained in the group. And that rather than it being proof that the `top-down' approach is somehow deeply significant it is rather a sign of the opposite: that the `bottom-up' approach constructs $H_4$ from $H_3$ whilst one can dualise the $30$ roots of $H_3$ directly to pure bivectors/quaternions using the pseudoscalar/inversion. This is not possible for $A_3$ where the inversion is missing and no pure quaternion representation of $A_3$ within $D_4$ exists; however, the Induction Theorem still holds and yields $A_3\rightarrow D_4$ \cite{Dechant2012CoxGA}.
\end{rem}

Having shown the existence and nature of various $H_4$ subrootsystems we now briefly discuss a nice way of visualising  4D polytopes  in the Coxeter plane, before using the subrootsystems of $H_4$ in order to construct pairs of invariant polytopes which we then visualise in the Coxeter plane. 

\section{The Coxeter Plane}\label{sec_CoxPl}

The Coxeter plane is a convenient way of visualising any root system in any dimension. The exposition is not necessary for the following sections but helps with the visualisation. We will briefly summarise the construction of this plane that is invariant under a corresponding Coxeter element. Its existence relies on the bipartite nature of the corresponding graphs (a two-colouring) \cite{Humphreys1990Coxeter}, which means that the simple roots can be partitioned into two mutually orthogonal sets (e.g. black and white), as can the reciprocal basis, the basis of fundamental weights. The properties of the Cartan matrix further mean that a Perron-Frobenius eigenvector with all positive entries exists. The components of this eigenvector corresponding to the black, respectively white, roots are used in a linear combination of the black, respectively white, fundamental weights. This gives a pair of (black and white) vectors which together determine a plane, which can be shown to be invariant under the Coxeter element. Since of course several such Coxeter elements exist (that are conjugate to one another), there are likewise several such planes. However, they give an equivalent description.  

The Clifford view of the Coxeter plane more generally has been investigated in \cite{Dechant2017e8}. Here, we instead perform all calculations in the 3D even subalgebra. The 4D simple roots can be chosen as follows in terms of the 3D $H_3$ simple roots:
\begin{align*}
\alpha_1&=a_1a_1&=1\\
\alpha_2&=a_1a_2&=-\frac{1}{2}  + \frac{1}{2}\tau e_{1}\wedge e_{2} + \frac{1}{2}\sigma e_{2}\wedge e_{3}\\
\alpha_3&=e_1e_2a_2e_3&=\frac{1}{2}\sigma e_{1}\wedge e_{2} -\frac{1}{2} e_{1}\wedge e_{3} + \frac{1}{2}\tau e_{2}\wedge e_{3}\\
\alpha_4&=a_2e_1e_2e_3&=\frac{1}{2}\sigma e_{1}\wedge e_{2} + \frac{1}{2} e_{1}\wedge e_{3}  -\frac{1}{2}\tau e_{2}\wedge e_{3}   \end{align*}

The reciprocal basis in this spinorial setup (with respect to the spinor inner product) is given by the following basis of fundamental weights:

\begin{align*}
   \omega_1&=1-\tau e_2\wedge e_3 - (\tau+1)e_1\wedge e_3\\
\omega_2&= -2\tau e_2\wedge e_3 - 2(\tau+1)e_1\wedge e_3\\
\omega_3&=-(2\tau+1)e_2\wedge e_3 - 3(\tau+1)e_1\wedge e_3 -\tau e_1\wedge e_2\\
\omega_4&=-(2\tau+1)e_2\wedge e_3 - (3\tau+1)e_1\wedge e_3 -\tau e_1\wedge e_2    \end{align*}

The Perron-Frobenius eigenvector of the $H_4$ Cartan matrix is given by $$v=\begin{pmatrix}
        4 + 4 \sqrt{5 } \\
        2 \left(1 + \sqrt{5}\right) \sqrt{\sqrt{5} + \sqrt{6 \sqrt{5} + 30} + 7}  \\
         \sqrt{6 \sqrt{5} + 30} + 8 + 4 \sqrt{5} + \sqrt{5} \sqrt{6 \sqrt{5} + 30}  \\
       \left(-1 + \sqrt{5} + \sqrt{6 \sqrt{5} + 30}\right) \sqrt{\sqrt{5} + \sqrt{6 \sqrt{5} + 30} + 7} \\
      \end{pmatrix}$$
 whilst the eigenvalue is
$$\lambda = 2- \frac{1}{2} \sqrt{7+\sqrt{5} + \sqrt{6 \sqrt{5} + 30}}. $$

\begin{center}
\begin{figure}
  
    \includegraphics[scale=0.75]{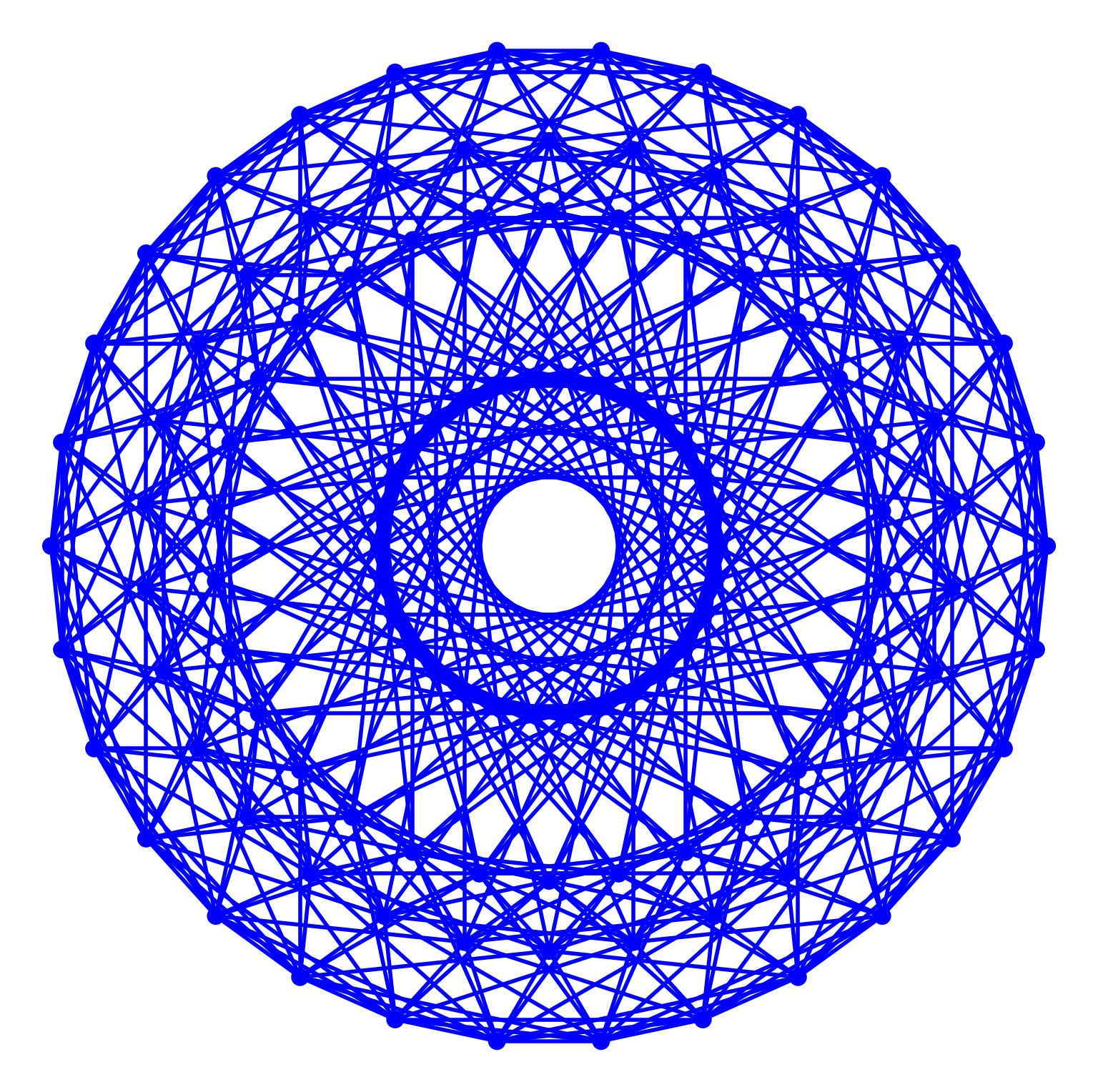}
    \caption{Projection of the $120$ roots of $H_4$ into the Coxeter plane. This is a familiar view -- however, all the calculations have been done in the even subalgebra of 3D. } 
    \label{fig:H4}
\end{figure}
\end{center}

Fig. \ref{fig:H4} shows the projection of the $120$ vertices and $720$ edges of the $H_4$ root system (aka the 600-cell) into its Coxeter plane. The projection is performed via the respective inner product with the `black and white spinors'. The computations are shown in the supplementary information.

We will use this Coxeter plane as a means to visualise the $H_4$ substructures of the following sections, including the Grand Antiprism and the snub 24-cell with their $H_2\oplus H_2$ and $D_4$ (aka $2T$) complements, as well as analogous constructions with $A_1\oplus A_1 \oplus A_1\oplus A_1$ ,  $A_2\oplus A_2$ and $A_4$ and their complements in the 600-cell.

\section{The Grand Antiprism and $H_2\times H_2$ }\label{sec_GA}

It is of course simple to show that $H_3$ contains an $H_2$ root system  (generated by the $a_2$ and $a_3$ simple roots), which leads to a corresponding root system $H_2$ in the 4D space of spinors. However, since the inversion $e_1e_2e_3$ is contained in the group $H_3$ this gets doubled to two orthogonal copies $H_2\oplus H_2$ sitting inside the $H_4$ root system, as seen above. A possible set is shown below, which is the one multiplicatively generated by $a_2$ and $a_3$  in combination with $e_1e_2e_3$ via the geometric product. Of course the $H_3$ root system, the icosidodecahedron, contains many such decagonal circles, but for this set of simple roots this set could be considered preferred:

\begin{center}
\begin{tabular}{|c|c|}
   \hline
\textbf{Element} $\times 2$ & \textbf{Element} $\times 2$  \\
  \hline 
  $ 2  $  & $ 2  e_{2}\wedge  e_{3}  $ \\
$ -2  $ & $ -2  e_{2}\wedge  e_{3}  $ \\
$ -\tau  -  e_{1}\wedge  e_{2} + \sigma  e_{1}\wedge  e_{3}  $ & $ \sigma  e_{1}\wedge  e_{2} +   e_{1}\wedge  e_{3} -\tau  e_{2}\wedge  e_{3}  $  \\
$ \sigma  -\tau  e_{1}\wedge  e_{2} -  e_{1}\wedge  e_{3}  $ & $ - e_{1}\wedge  e_{2} + \tau  e_{1}\wedge  e_{3} + \sigma  e_{2}\wedge  e_{3}  $  \\
$ \tau  -  e_{1}\wedge  e_{2} + \sigma  e_{1}\wedge  e_{3}  $ & $ \sigma  e_{1}\wedge  e_{2} +   e_{1}\wedge  e_{3} + \tau  e_{2}\wedge  e_{3}  $  \\
$ \tau  +   e_{1}\wedge  e_{2} -\sigma  e_{1}\wedge  e_{3}  $ & $ -\sigma  e_{1}\wedge  e_{2} -  e_{1}\wedge  e_{3} + \tau  e_{2}\wedge  e_{3}  $  \\
$ -\sigma  + \tau  e_{1}\wedge  e_{2} +   e_{1}\wedge  e_{3}  $ & $   e_{1}\wedge  e_{2} -\tau  e_{1}\wedge  e_{3} -\sigma  e_{2}\wedge  e_{3}  $ \\
$ -\sigma  -\tau  e_{1}\wedge  e_{2} -  e_{1}\wedge  e_{3}  $ & $ -  e_{1}\wedge  e_{2} + \tau  e_{1}\wedge  e_{3} -\sigma  e_{2}\wedge  e_{3}  $  \\
$ -\tau  +   e_{1}\wedge  e_{2} -\sigma  e_{1}\wedge  e_{3}  $ & $ -\sigma  e_{1}\wedge  e_{2} -  e_{1}\wedge  e_{3} -\tau  e_{2}\wedge  e_{3}  $  \\
$ \sigma  + \tau  e_{1}\wedge  e_{2} +   e_{1}\wedge  e_{3}  $ &$   e_{1}\wedge  e_{2} -\tau  e_{1}\wedge  e_{3} + \sigma  e_{2}\wedge  e_{3}  $ \\

  \hline 

\end{tabular}
\end{center}

  \begin{center}
\begin{figure}
  
    \includegraphics[scale=0.75]{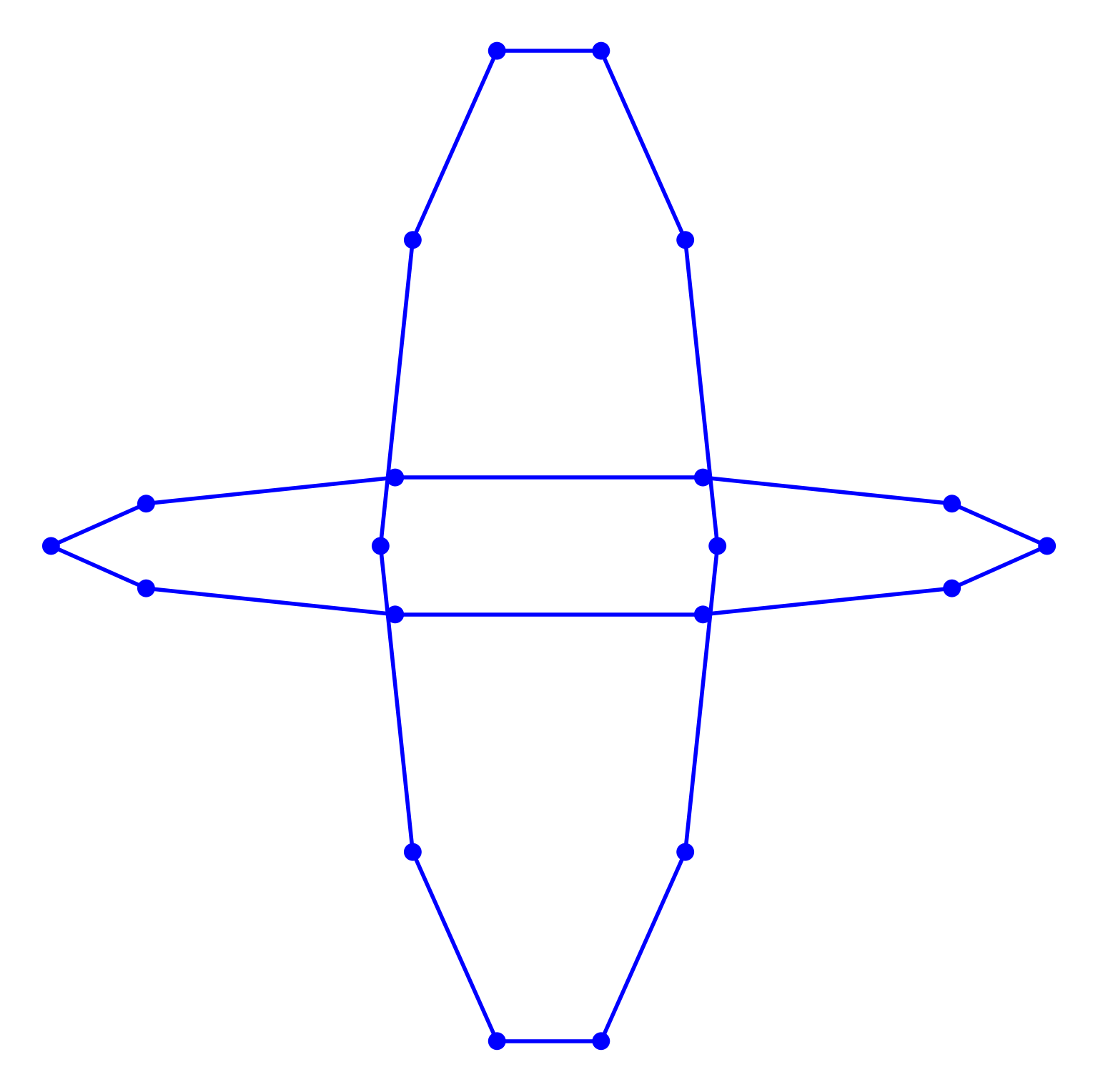}
    \caption{The projection of the $20$ roots of the $H_2\oplus H_2$ sitting inside the 600-cell/$H_4$ root system into the $H_4$ Coxeter plane.}   \label{fig_H2H2}
\end{figure}  
\end{center}

Of course they are invariant under their own automorphism group $\Aut(H_2\oplus H_2)$ of order $20\times 20 = 400$. In previous work the author has already argued that the automorphism group of a spinorial/quaternionic root system is just the group acting on itself by left and right multiplication, leading to two factors of the same group \cite{Dechant2013Platonic}. This spinorial group multiplication is to be distinguished from the above `spin reflections', which one could consider a different type of multiplication, generating the reflection/Coxeter groups, which was termed `conjugal' in the above paper. The $H_2\oplus H_2$ root system is invariant under this conjugal group multiplication (i.e. 4D reflections) by virtue of being a spinor group by Proposition \ref{prop_spin_closure}.

  \begin{center}
\begin{figure}
  
    \includegraphics[scale=0.75]{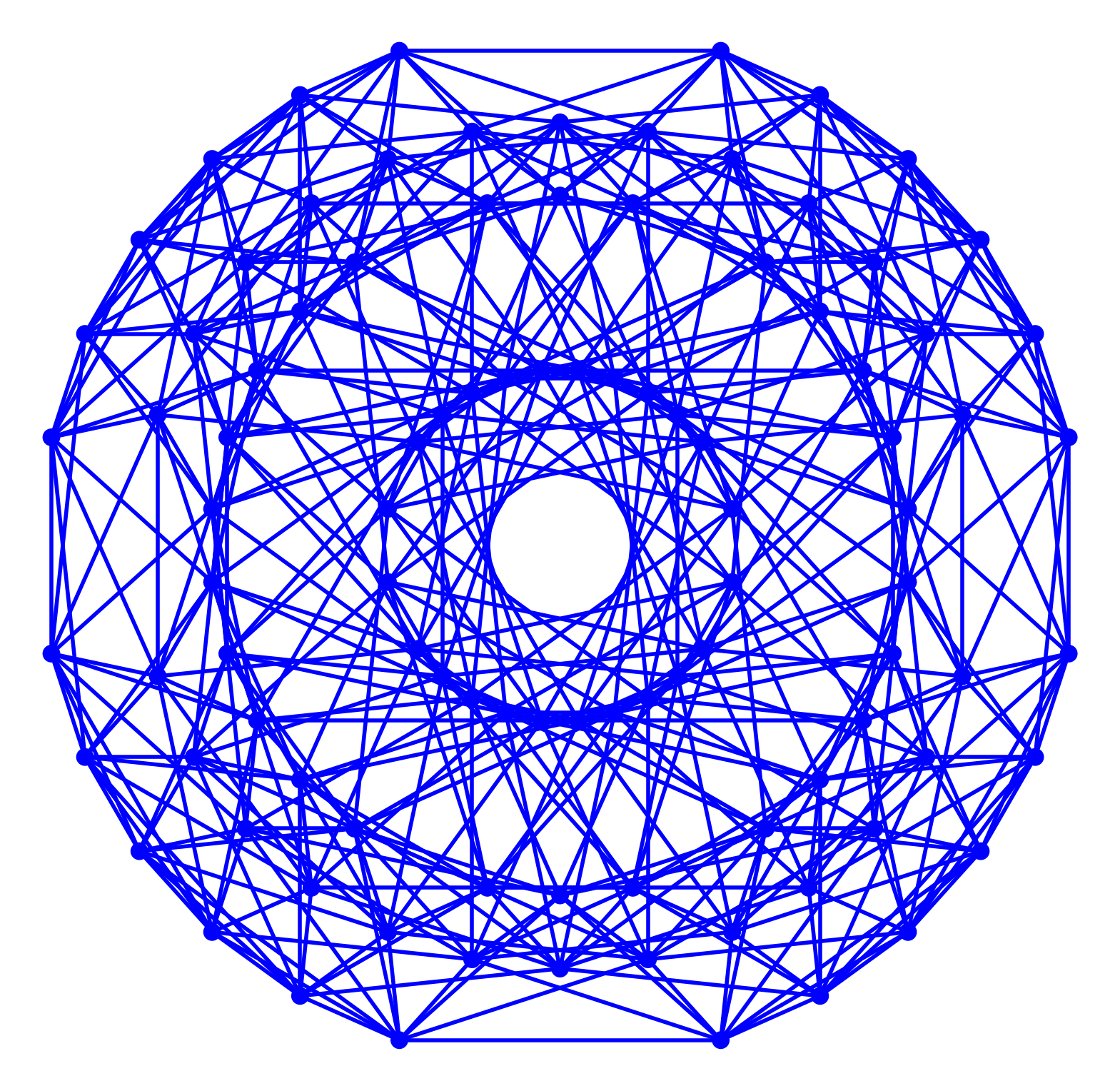}
  \caption{Removal of those $20$ roots of the $H_2\oplus H_2$ from the 600-cell/$H_4$ root system leads to the Grand Antiprism with 100 vertices. Their projection into the $H_4$ Coxeter plane is shown above.}     \label{fig_GA}
\end{figure}  \end{center}

But since this $H_2\times H_2$ is a subgroup of $H_4$,  its complement in the 600-cell is separately left invariant. By this we mean the collection of vertices derived from the 600-cell by subtracting from the $120$ vertices the $20$ vertices from the $H_2\oplus H_2$ root system. This orbit / collection of points is separately invariant and therefore has the same automorphism group of order $400$. This can also be verified straightforwardly by direct computation (see supplementary material). This complement of the $H_2\oplus H_2$ root system is a 4D polytope with $100$ vertices and $500$ edges which was found in 1965 in \cite{Conway1965Four} by computational means and is called the `Grand Antiprism'. This is also discussed from a quaternionic perspective in \cite{Koca2009grand, Koca2006H4,koca2007group,koca2012snub}. The construction from 3D seems both simpler in terms of deriving the vertex set and in shedding light on its symmetry group, as well as the conceptual and uniform construction that carries over in the other cases below. The projection of the two orthogonal $H_2$s into the $H_4$ Coxeter plane is shown in Fig. \ref{fig_H2H2}. The corresponding projection of the Grand Antiprism is shown in Fig. \ref{fig_GA}. For ease of visualisation the edges are also computed and projected, and plotted in 3D in a SageMath visualisation that can be further explored \cite{stein2008sage} (supplementary information).

\begin{prop}[$H_2\times H_2$ split of $H_4$ vertices]
The sets $H_2\oplus H_2$ and its complement the Grand Antiprism are separately invariant under $H_2\times H_2$. \end{prop}

\begin{proof}
Straightforward or by straightforward explicit calculation (supplementary material). 
\end{proof}

\section{Snub 24-cell }\label{sec_snub}

\begin{center}
\begin{figure}
  
    \includegraphics[scale=0.75]{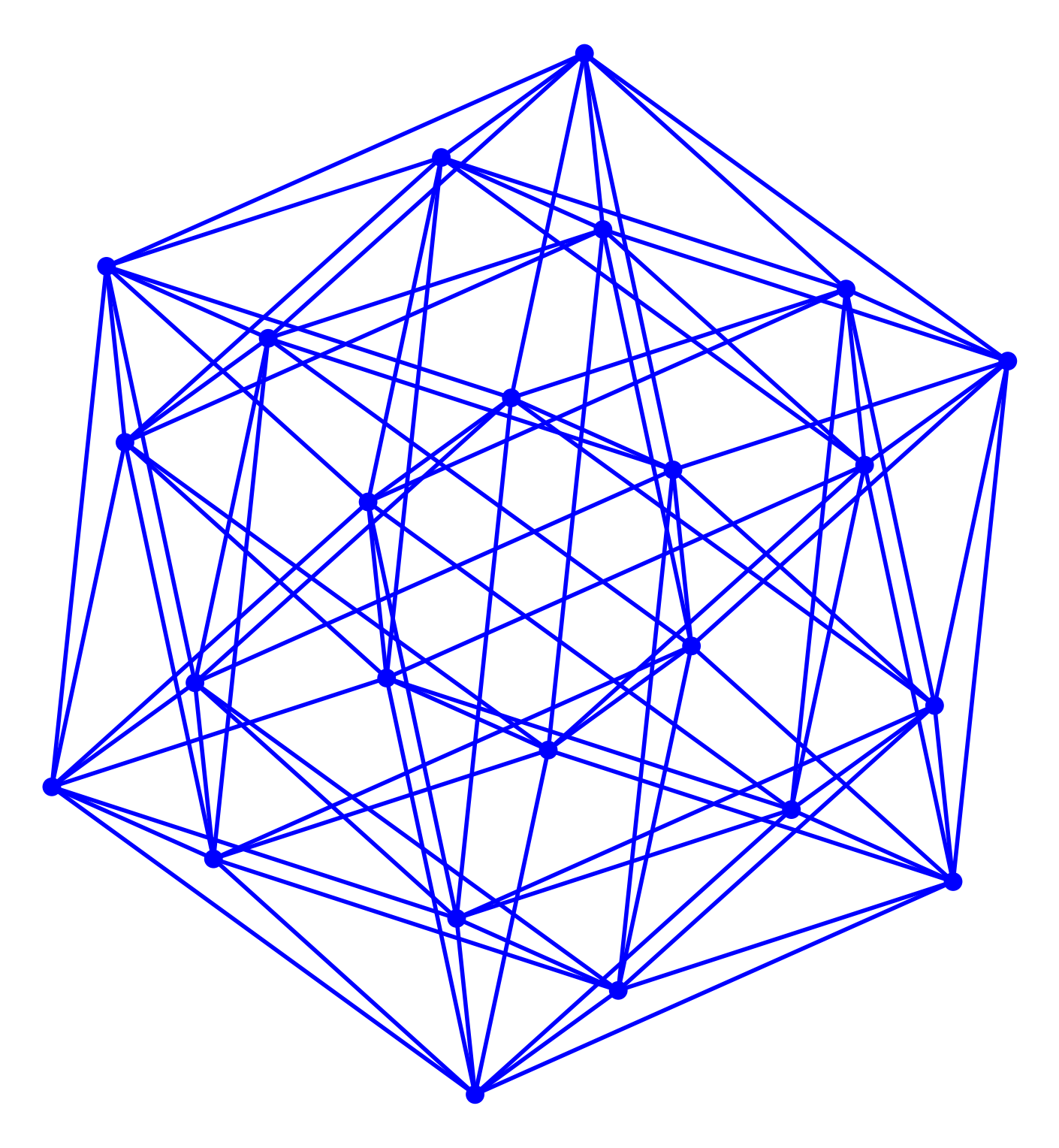}
    \caption{The projection of the $24$ roots of the $D_4$ root system sitting inside the 600-cell/$H_4$ root system into the $H_4$ Coxeter plane. These vertices of course come from the binary tetrahedral group $2T$ inside the $2I$.} 
    \label{fig:2T}
\end{figure}
\end{center}

\begin{center}
\begin{figure}

    \includegraphics[scale=0.75]{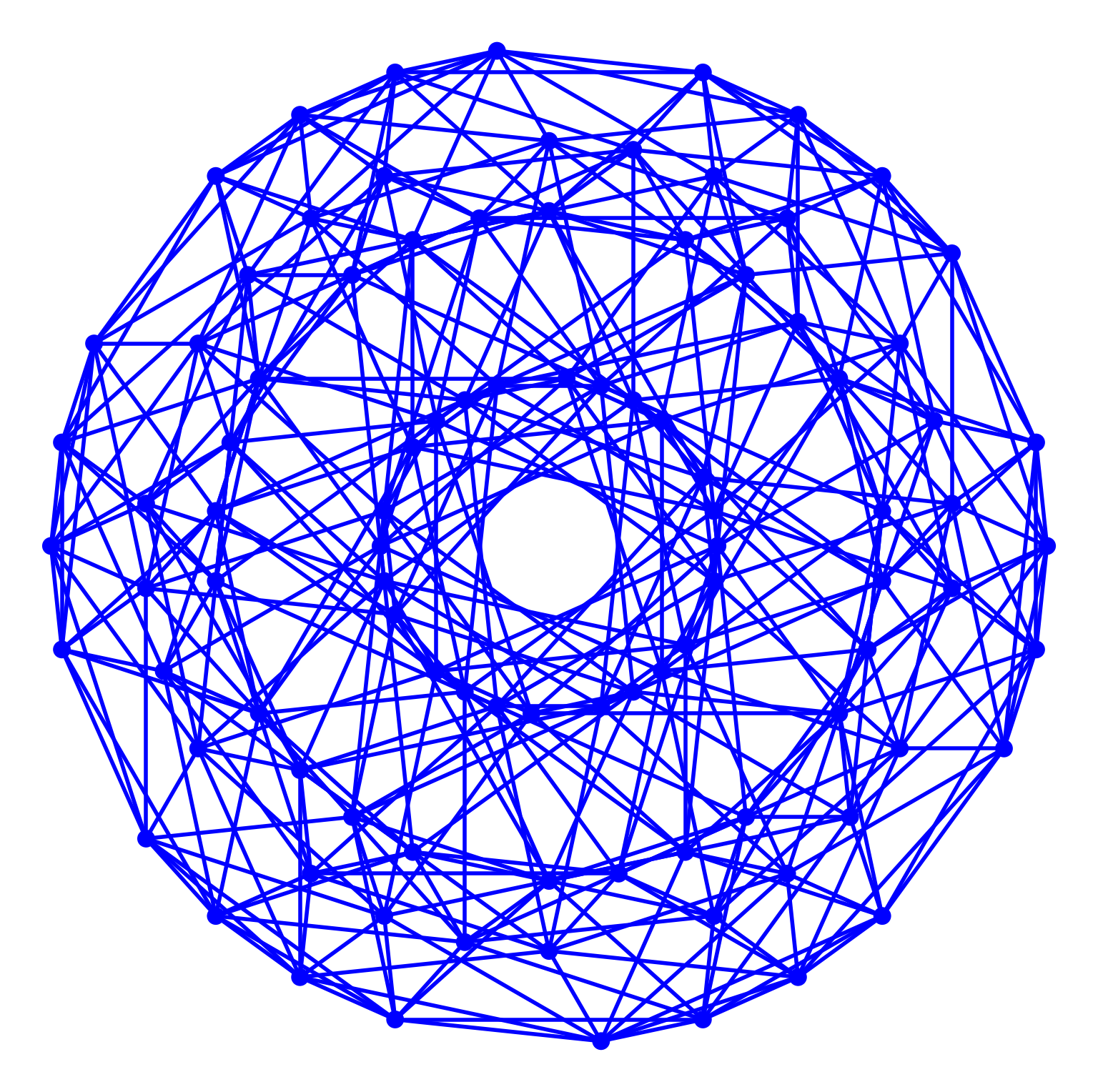}
    \caption{The projection of the snub 24-cell into the $H_4$ Coxeter plane. The 96 vertices can be derived by removing the $24$ vertices of $D_4$/the binary tetrahedral group from the 600-cell.} 
    \label{fig:snub}
\end{figure}
\end{center}

From a 3D perspective it is obvious that although $A_3$ is not a subrootsystem of $H_3$, the binary tetrahedral group $2T$ \emph{is} contained in the binary icosahedral group $2I$ (this can also easily been seen since the tetrahedral group is the alternating group $A_4$ whilst the icosahedral group is $A_5$ and thus the former is contained in the latter, and this also holds for their spin double covers) and the root system $D_4$ is thus contained in $H_4$. As a spin/quaternionic group, $D_4$ is of course just $2T$ and its automorphism group is just $2T\times 2T$ of order $24^2=576$. Removing the $24$ vertices from the 600-cell leads to a set of $96$ points that is separately invariant under $D_4$, in analogy to the construction of the Grand Antiprism above. This collection of $96$ vertices connected by $432$ edges is known as the `snub 24-cell'. Again we believe the construction from 3D to be conceptually clearer, more systematic and more economical. The projection of the $24$ $2T$ spinors (listed below) aka the $D_4$ root system  along with its  $96$ edges into the $H_4$ Coxeter plane is shown in Fig. \ref{fig:2T}. The corresponding projection of the snub 24-cell and its edges is shown in Fig. \ref{fig:snub}.

\begin{center}
\begin{tabular}{|c|c|}
   \hline
\textbf{Element} $\times 2$ & \textbf{Element} $\times 2$  \\
  \hline 
$2$& $-2$\\
 $e_1\wedge e_2 -\tau e_1\wedge e_3 + \sigma e_2\wedge e_3$& $-e_1\wedge e_2 + \tau e_1\wedge e_3 -\sigma e_2\wedge e_3$\\

 $\sigma e_1\wedge e_2 - e_1\wedge e_3 + \tau e_2\wedge e_3$&$-\sigma e_1\wedge e_2 + e_1\wedge e_3 -\tau e_2\wedge e_3$\\

 $\tau e_1\wedge e_2 -\sigma e_1\wedge e_3 + e_2\wedge e_3$& $-\tau e_1\wedge e_2 + \sigma e_1\wedge e_3 - e_2\wedge e_3$\\

$1 + e_1\wedge e_2 - e_1\wedge e_3 + e_2\wedge e_3$&  $-1 - e_1\wedge e_2 + e_1\wedge e_3 - e_2\wedge e_3$\\
  $1 - e_1\wedge e_2 + e_1\wedge e_3 - e_2\wedge e_3$& $-1 + e_1\wedge e_2 - e_1\wedge e_3 + e_2\wedge e_3$\\
 $1 + \tau e_1\wedge e_2 + \sigma e_2\wedge e_3$&$-1 -\tau e_1\wedge e_2 -\sigma e_2\wedge e_3$\\
$1 -\tau e_1\wedge e_2 -\sigma e_2\wedge e_3$&$-1 + \tau e_1\wedge e_2 + \sigma e_2\wedge e_3$\\

 $1 + \sigma e_1\wedge e_2 -\tau e_1\wedge e_3$& $-1 -\sigma e_1\wedge e_2 + \tau e_1\wedge e_3$\\
$1 -\sigma e_1\wedge e_2 + \tau e_1\wedge e_3$& $-1 + \sigma e_1\wedge e_2 -\tau e_1\wedge e_3$\\
 
 $1 + \sigma e_1\wedge e_3 -\tau e_2\wedge e_3$&$-1 -\sigma e_1\wedge e_3 + \tau e_2\wedge e_3$\\
  $1 -\sigma e_1\wedge e_3 + \tau e_2\wedge e_3$& $-1 + \sigma e_1\wedge e_3 -\tau e_2\wedge e_3$\\

\hline 
\end{tabular}
\end{center}

Note that because of the spin double cover property they come at least in pairs. But if they are not their own reverse,  then they even come in quadruplets related via reversal and multiplication by $-1$ (c.f. Proposition \ref{prop_spin_closure}).

\begin{prop}[$D_4$ split of $H_4$ vertices]
The sets $D_4$ and its complement the snub 24-cell are separately invariant under $D_4$. \end{prop}

\begin{proof}
Straightforward or by straightforward explicit calculation (supplementary material). 
\end{proof}

\section{$A_1^4$, $A_2\times A_2$ and $A_4$}\label{sec_A2}

We follow the above uniform construction of splitting the $H_4$ root system with respect to its subrootsystems, using the remaining examples of $A_1\oplus A_1\oplus A_1\oplus A_1$, $A_2\oplus A_2$ and $A_4$.

  \begin{center}
\begin{figure}
  
    \includegraphics[scale=0.75]{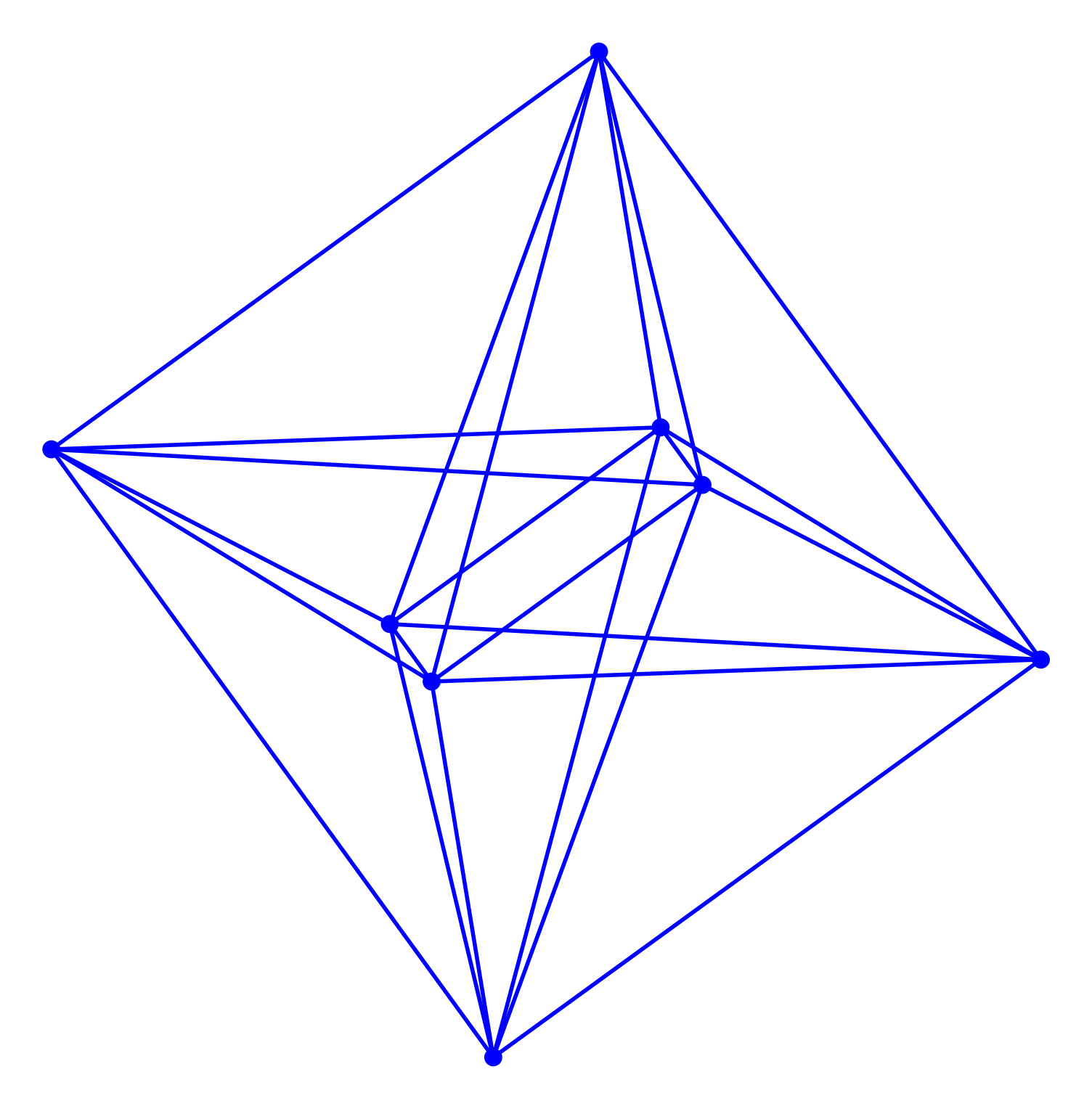}
    \caption{The projection of the $8$ roots of the $A_1^4$ sitting inside the 600-cell/$H_4$ root system into the $H_4$ Coxeter plane.} 

    \label{fig:A14}
\end{figure}
\end{center}

\begin{center}
\begin{figure}
  
    \includegraphics[scale=0.75]{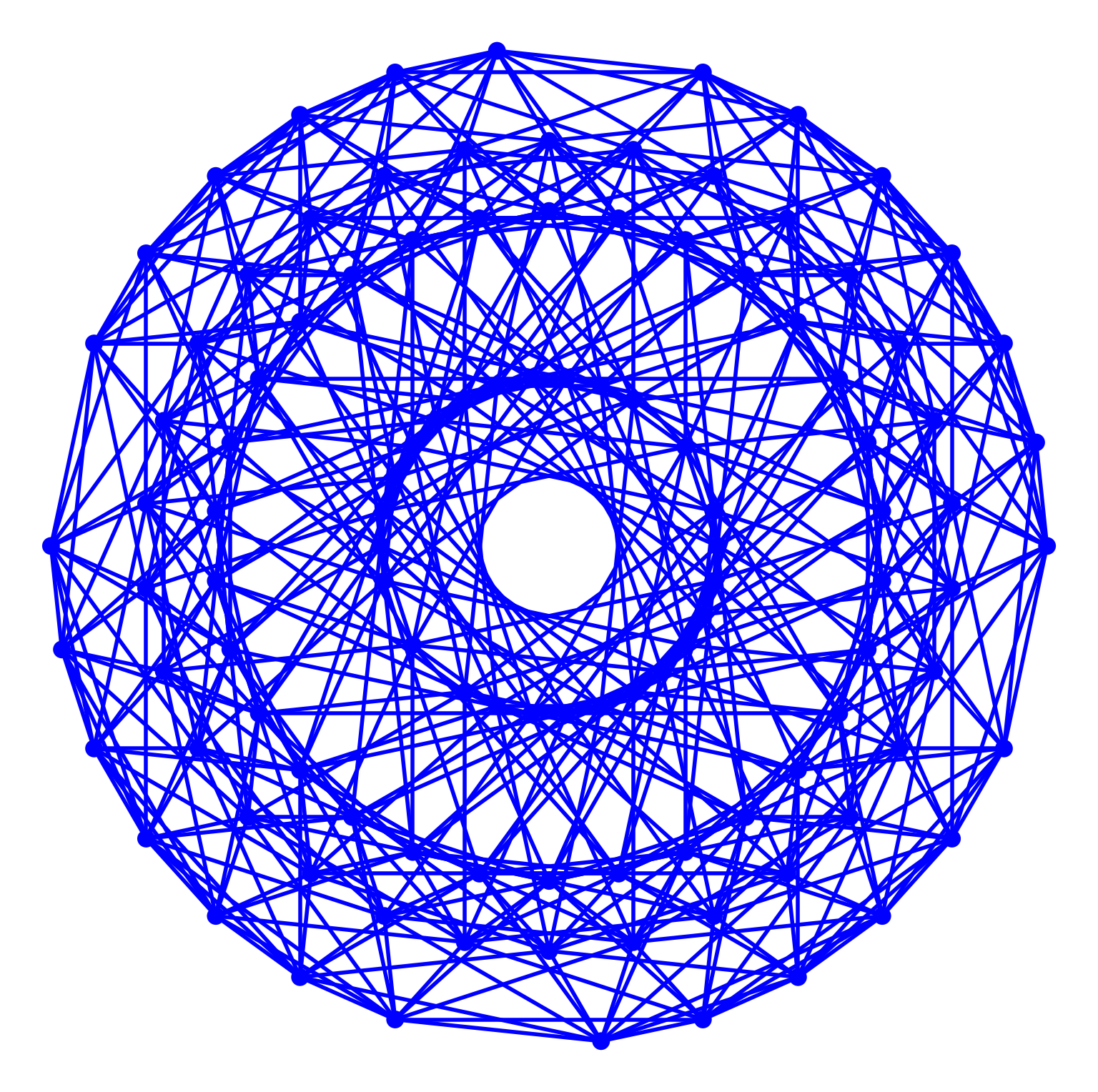}
    \caption{The projection of the $112$ vertices into the $H_4$ Coxeter plane that are left from the 600-cell by removing the above $8$ vertices constituting the $A_1^4$  root system, along with its 624 edges.} 
    \label{fig_CA14}
\end{figure}
\end{center}

The basic unit vectors $e_1$, $e_2$, $e_3$ of course generate the quaternion group consisting of $\pm 1$, $\pm e_1 e_2$, $\pm e_2 e_3$ and $\pm e_3 e_1$. The projection of this 16-cell with its $8$ vertices and $24$ edges into the $H_4$ Coxeter plane is shown in Fig. \ref{fig:A14}. Of course its complement is also invariant under $A_1\times A_1\times A_1\times A_1$:

\begin{prop}[$A_1\times A_1\times A_1\times A_1$ split of $H_4$ vertices]
The sets  $A_1\oplus A_1\oplus A_1\oplus A_1$ and its complement are separately invariant under $A_1\times A_1\times A_1\times A_1$. \end{prop}
\begin{proof}
Straightforward or by straightforward explicit calculation (supplementary material). 
\end{proof}
The projection of this complement of this 16-cell with its $112$ vertices and $624$ edges into the $H_4$ Coxeter plane is shown in Fig. \ref{fig_CA14}.

Similarly to the $H_2$ case above, it is simple to show that $H_3$ contains an $A_2$ root system (generated via reflections in the $a_1$ and $a_2$ simple roots), which leads to a corresponding root system $A_2$ in the 4D space of spinors. However, since the inversion $e_1e_2e_3$ is contained in the group $H_3$ this $A_2$ again  gets doubled to two orthogonal copies $A_2\oplus A_2$ sitting inside the $H_4$ root system, as seen above. A possible set is shown below, which is the one multiplicatively generated by $a_1$ and $a_2$ simple roots of $H_3$ in combination with the inversion $e_1e_2e_3$. But other choices would be possible, since of course the $H_3$ root system, the icosidodecahedron, contains many such hexagonal circles.

\begin{center}
\begin{tabular}{|c|c|}
   \hline
\textbf{Element} $\times 2$ & \textbf{Element} $\times 2$  \\
  \hline 
 $ 2  $ & $ -2  e_{1}\wedge  e_{3}  $  \\
$ -2  $ & $ 2  e_{1}\wedge  e_{3}  $  \\
$ -1  + \tau  e_{1}\wedge  e_{2} + \sigma  e_{2}\wedge  e_{3}  $ & $ \sigma  e_{1}\wedge  e_{2} + 1  e_{1}\wedge  e_{3} -\tau  e_{2}\wedge  e_{3}  $ \\
$ 1  + \tau  e_{1}\wedge  e_{2} + \sigma  e_{2}\wedge  e_{3}  $ & $ \sigma  e_{1}\wedge  e_{2} -1  e_{1}\wedge  e_{3} -\tau  e_{2}\wedge  e_{3}  $ \\
$ 1  -\tau  e_{1}\wedge  e_{2} -\sigma  e_{2}\wedge  e_{3}  $ & $ -\sigma  e_{1}\wedge  e_{2} -1  e_{1}\wedge  e_{3} + \tau  e_{2}\wedge  e_{3}  $ \\
$ -1  -\tau  e_{1}\wedge  e_{2} -\sigma  e_{2}\wedge  e_{3}  $ & $ -\sigma  e_{1}\wedge  e_{2} + 1  e_{1}\wedge  e_{3} + \tau  e_{2}\wedge  e_{3}  $ \\

\hline 
\end{tabular}
\end{center}

  \begin{center}
\begin{figure}
  
    \includegraphics[scale=0.75]{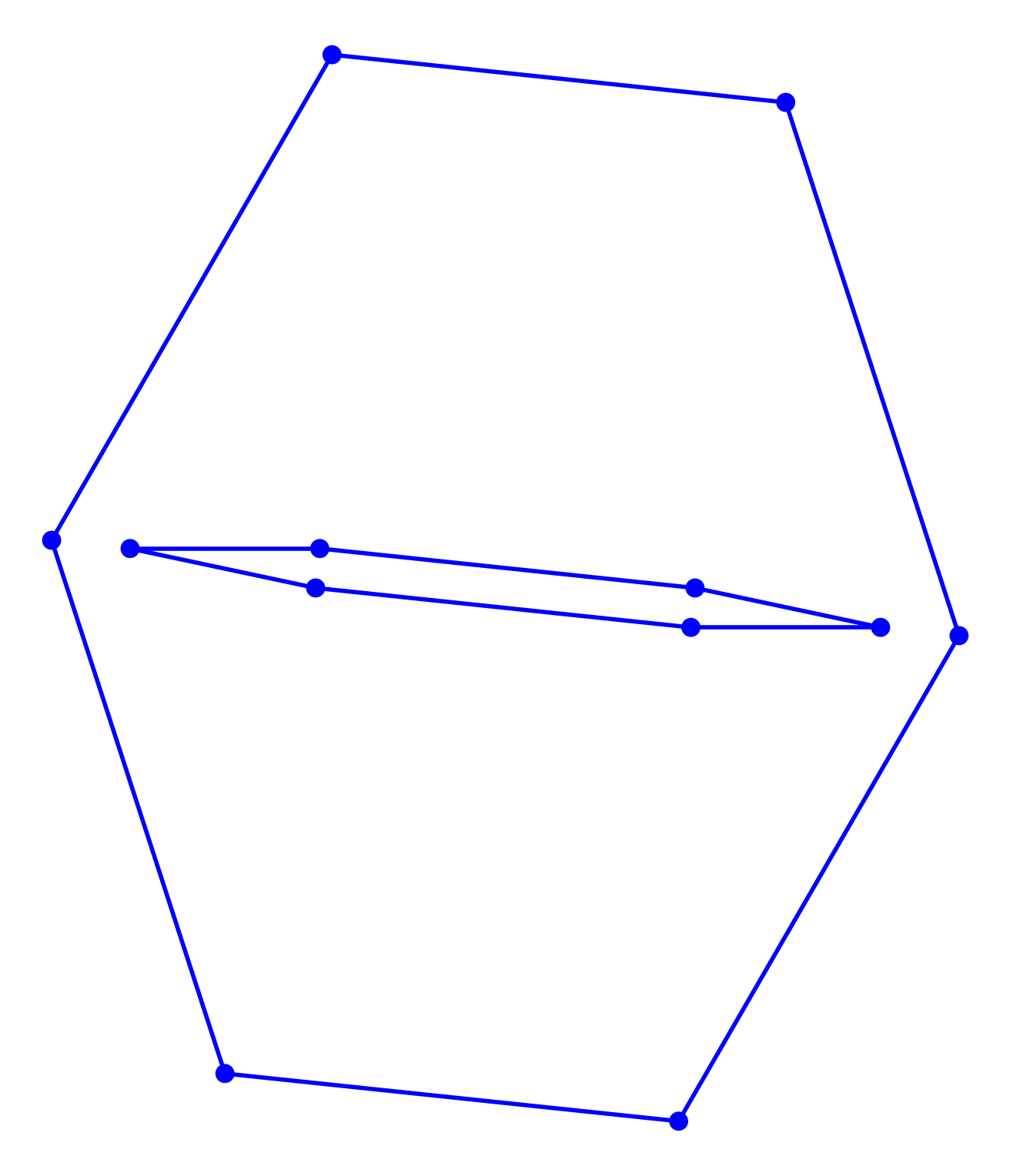}
    \caption{The projection of the $12$ roots of the $A_2\oplus A_2$ sitting inside the 600-cell/$H_4$ root system into the $H_4$ Coxeter plane.} 

    \label{fig:A2A2}
\end{figure}
\end{center}

\begin{center}
\begin{figure}
  
    \includegraphics[scale=0.75]{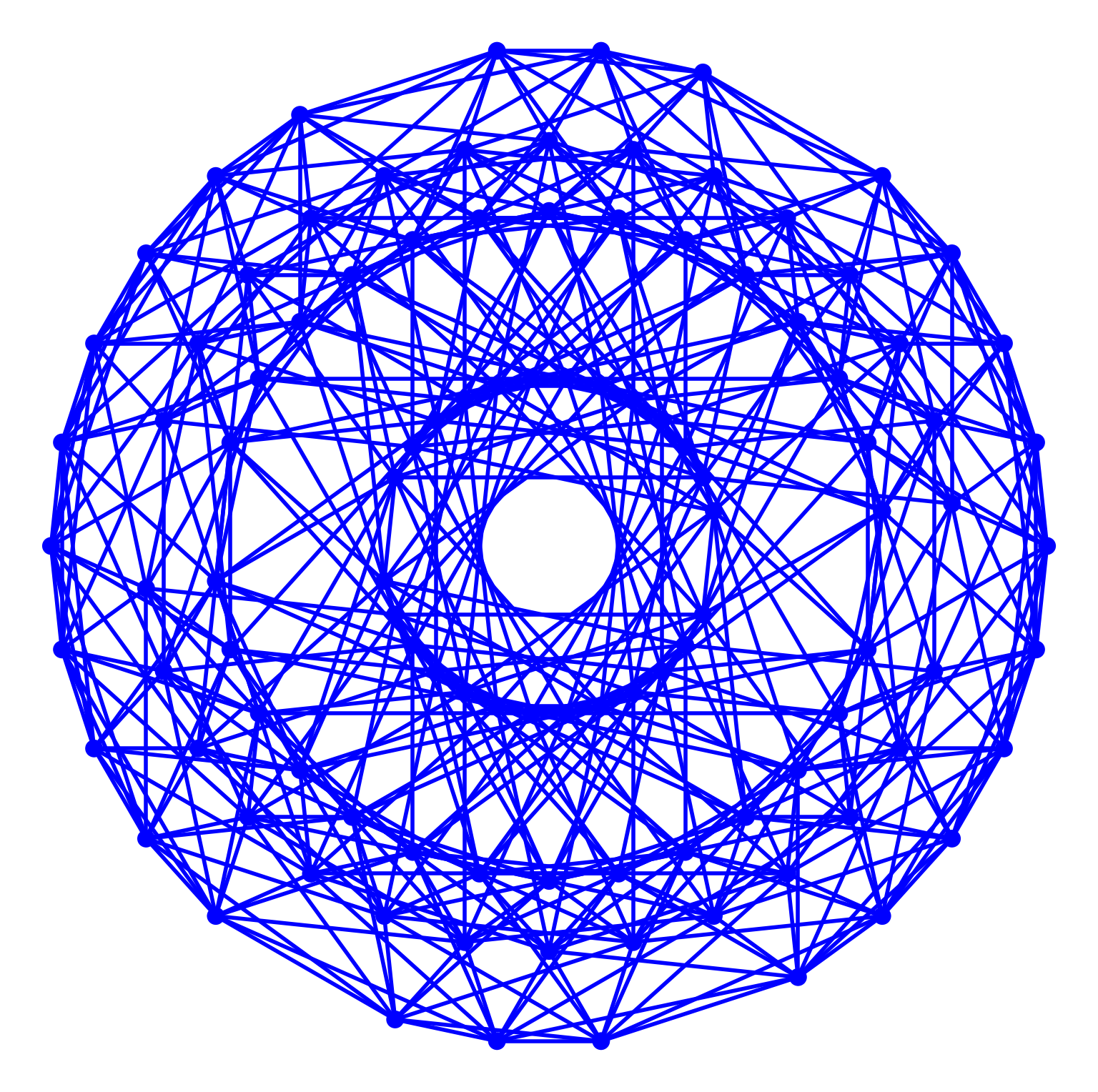}
    \caption{The projection of the $108$ vertices into the $H_4$ Coxeter plane that are left from the 600-cell by removing the above $12$ vertices from the $A_2\oplus A_2$  root system.} 
    \label{fig:CA2A2}
\end{figure}
\end{center}

The projection of these $12$ points into the $H_4$ Coxeter plane as before is shown in Fig. \ref{fig:A2A2}. Its complement consisting of $108$ vertices and $576$ edges is shown in Fig. \ref{fig:CA2A2}.

\begin{prop}[$A_2\times A_2$ split of $H_4$ vertices]
The sets $A_2\oplus A_2$ and its complement are separately invariant under $A_2\times A_2$. \end{prop}

\begin{proof}
Straightforward or by straightforward explicit calculation (supplementary material). 
\end{proof}

Our final example is the $A_4$ root system contained in $H_4$ as we saw above. Again the $A_4$ root system and its complement are both invariant.

\begin{prop}[$A_4$ split of $H_4$ vertices]
The sets $A_4$ and its complement are separately invariant under $A_4$. \end{prop}

\begin{proof}
Straightforward or by straightforward explicit calculation (supplementary material). 
\end{proof}

\begin{center}
\begin{figure}
  
    \includegraphics[scale=0.75]{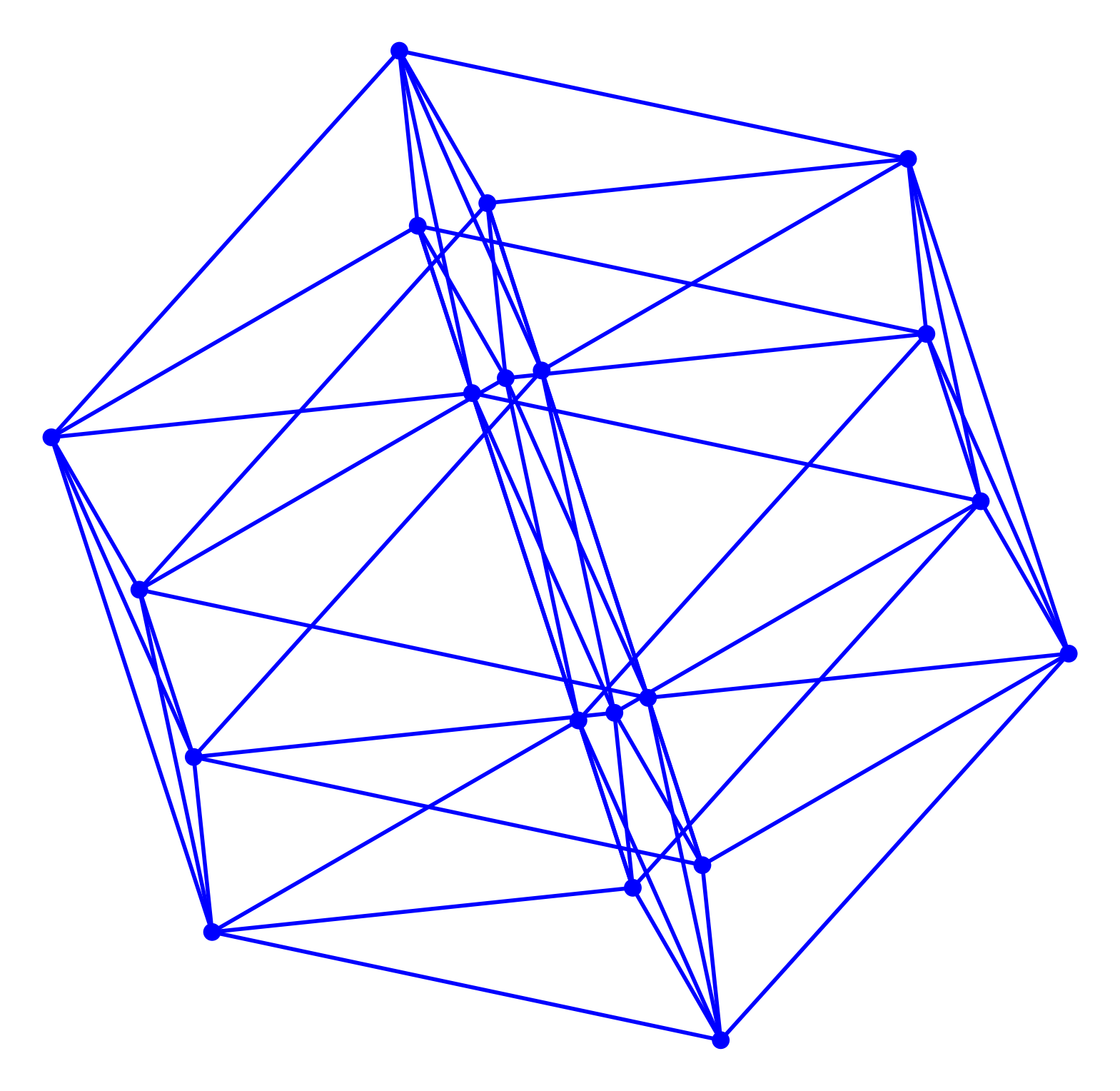}
    \caption{The projection of the $20$ roots of the $A_4$ root system sitting inside the 600-cell/$H_4$ root system into the $H_4$ Coxeter plane. } 
    \label{fig:A4}
\end{figure}
\end{center}

\begin{center}
\begin{figure}

    \includegraphics[scale=0.75]{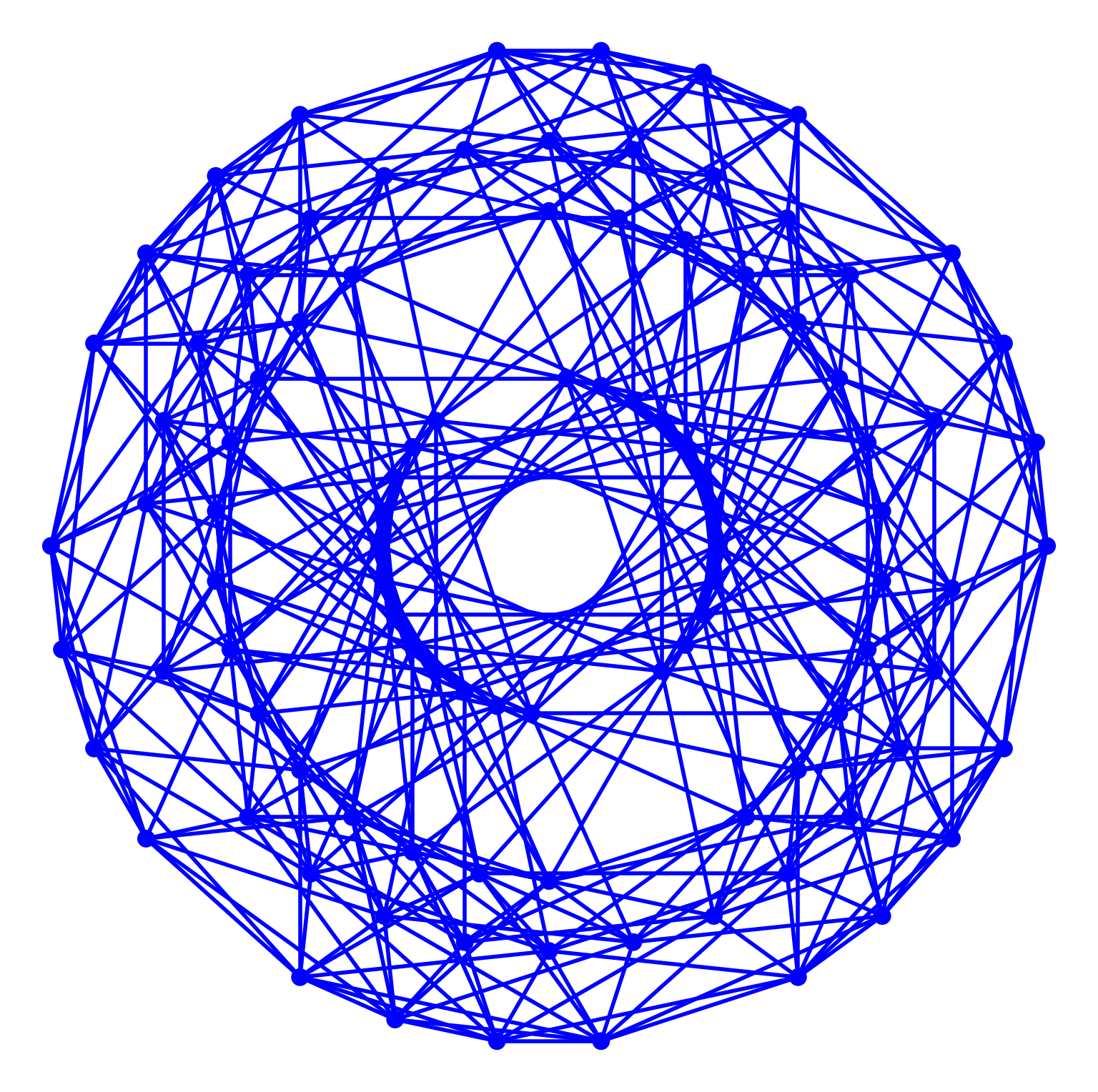}
    \caption{The projection of the complement of the above into the $H_4$ Coxeter plane. The 100 vertices can be derived by removing the $20$ vertices of $A_4$ from the 600-cell.} 
    \label{fig:CA4}
\end{figure}
\end{center}

The root system consists of $20$ vertices and $60$ edges (shown in Fig. \ref{fig:A4}) whilst its complement consists of $100$ vertices and $480$ edges and is shown in  Fig. \ref{fig:CA4}.

A possible choice of $A_4$ roots is shown here:
\begin{center}
\begin{tabular}{|c|c|}
   \hline
\textbf{Element} $\times 2$ & \textbf{Element} $\times 2$  \\
  \hline 
 $2$&$-2$\\
 
$\tau e_1\wedge e_2 + \sigma e_1\wedge e_3 + e_2\wedge e_3$& $-\tau e_1\wedge e_2 -\sigma e_1\wedge e_3 - e_2\wedge e_3$\\
 $e_1\wedge e_2 -\tau e_1\wedge e_3 + \sigma e_2\wedge e_3$&$-e_1\wedge e_2 + \tau e_1\wedge e_3 -\sigma e_2\wedge e_3$\\
 $\sigma e_1\wedge e_2 - e_1\wedge e_3 -\tau e_2\wedge e_3$& $-\sigma e_1\wedge e_2 + e_1\wedge e_3 + \tau e_2\wedge e_3$\\

 $1 + \tau e_1\wedge e_2 + \sigma e_2\wedge e_3$&$-1 -\tau e_1\wedge e_2 -\sigma e_2\wedge e_3$\\
 $1 -\tau e_1\wedge e_2 -\sigma e_2\wedge e_3$&$-1 + \tau e_1\wedge e_2 + \sigma e_2\wedge e_3$\\

  $1 + \sigma e_1\wedge e_3 + \tau e_2\wedge e_3$& $-1 -\sigma e_1\wedge e_3 -\tau e_2\wedge e_3$\\
 $1 -\sigma e_1\wedge e_3 -\tau e_2\wedge e_3$& $-1 + \sigma e_1\wedge e_3 + \tau e_2\wedge e_3$\\
 
 $1 + \sigma e_1\wedge e_2 -\tau e_1\wedge e_3$& $-1 -\sigma e_1\wedge e_2 + \tau e_1\wedge e_3$\\
$1 -\sigma e_1\wedge e_2 + \tau e_1\wedge e_3$& $-1 + \sigma e_1\wedge e_2 -\tau e_1\wedge e_3$\\
 
\hline 
\end{tabular}
\end{center}

Note that because of the spin double cover property these elements come at least in pairs (which satisfies the corresponding root system axiom). But if they are not their own reverse,  then they even come in quadruplets related via reversal \emph{and} multiplication by $-1$ (c.f. Proposition \ref{prop_spin_closure}).  

This concludes our listed examples, which illustrate the uniform approach of constructing subrootsystems via the connection with 3D and splitting the vertices of the 600-cell into two separately invariant sets. The $H_4$ Coxeter plane provides a nice visualisation for each complementary pair based on $H_4$ subrootsystems. 

\section{Conclusions}\label{sec_concl}

The intention for this article was to firstly shed light on 4D geometry, root systems and polytopes through the connection with 3D spinors via the uniform Induction Theorem. This gives additional insight into 4D root systems and polytopes along with their symmetries via another uniform construction that splits the $H_4$ root system into a complementary pair of separately invariant polytopes, consisting of a subrootsystem and its complement in $H_4$. These can also be consistently visualised via a projection into the Coxeter plane. As we have shown in previous work there are many connections across exceptional objects throughout mathematics, including Trinities and ADE correspondences, which this work relates to. In particular, the `bottom-up' view of exceptional objects has led to a very fruitful and insightful line of research, and is perhaps mirrored by some constructions in finite group theory for some of the sporadic groups. We continue to advocate the use of geometric insight in addition to the purely algebraic manipulation in terms of quaternions, by viewing quaternions as arising in a geometric guise as spinors in 3 dimensions with a much clearer and consistent geometric interpretation.   

Secondly, we continue to advocate that Clifford algebras and root systems/reflection groups are  natural and complementary frameworks - in a setting with a vector space with an inner product, and a powerful reflection formula - and can therefore be synthesised into one powerful and coherent framework. Thirdly, this provides a general arena to do group theoretical calculations in as demonstrated with some detailed examples in this paper, and we will investigate group and representation theoretic topics in more detail in future work. In particular, in the spirit of open science and reproducibility I hope that the computational work sheets provided as supplementary material are useful and can help further research in this area by adapting my code and using other free software provided by the community such as \verb|galgebra| and Sage. Finally, at the end of 2020, we wish to commemorate and honour some of the giants in this field, Conway and the Guys. 


\subsection*{Acknowledgment}
This paper is dedicated to the memory of the late John Horton Conway and Richard Kenneth Guy, as well as Michael Guy.  I would like to thank Yang-Hui He, Yvan Saint-Aubin, Ji\v ri Patera, Peter Cameron, John McKay, Robert Wilson, Alexander Konovalov, Jean-Jacques Dupas, Tony Sudbery, Reidun Twarock, David Hestenes, Anthony Lasenby, Joan Lasenby and Eckhard Hitzer for interesting discussions over the years, and Hugo Hadfield and Eric Wieser for help with setting up the \verb|galgebra| software package \cite{Bromborsky2020}. 


\bibliography{virobib}

\end{document}